\documentclass[12pt]{article}  
\usepackage{amsthm,amsmath,amsfonts,epsfig,amssymb,latexsym}
\usepackage{mathrsfs}
\usepackage[T1]{fontenc}
\usepackage{graphicx}
\usepackage{enumerate}
\usepackage{xy,xypic}
\usepackage{graphics}
\usepackage{footnote}
\usepackage{scalefnt}
\usepackage{tikz}
\usetikzlibrary{shapes,arrows,fit,calc,positioning}
\usetikzlibrary{matrix}
\xyoption{all}
\title{{\large Convergence of two obstructions for projective modules}} 
\author{
Satya Mandal\footnote{Partially supported by a General Research Grant (no 2301857) from U. of Kansas}
 \\ 
{\small University of Kansas, Lawrence, Kansas 66045, USA}\\
{\small {\it  mandal@ku.edu} 
  }\\
 } 
\begin{document}
\renewcommand{\baselinestretch}{1.255}
\setlength{\parskip}{1ex plus0.5ex}
\date{2 October 2023}

\newtheorem{theorem}{Theorem}[section]
\newtheorem{proposition}[theorem]{Proposition}
\newtheorem{lemma}[theorem]{Lemma}
\newtheorem{corollary}[theorem]{Corollary}
\newtheorem{construction}[theorem]{Construction}
\newtheorem{notations}[theorem]{Notations}
\newtheorem{question}[theorem]{Question}
\newtheorem{example}[theorem]{Example}
\newtheorem{definition}[theorem]{Definition} 
\newtheorem{conjecture}[theorem]{Conjecture} 
\newtheorem{remark}[theorem]{Remark} 
\newtheorem{statement}[theorem]{Statement}

\newcommand{\iso}{\stackrel{\sim}{\longrightarrow}}
\def\spec#1{\mathrm{Spec}\left(#1\right)}

\newcommand{\sur}{\twoheadrightarrow}
\newcommand{\bD}{\begin{definition}}
\newcommand{\eD}{\end{definition}}
\newcommand{\bP}{\begin{proposition}}
\newcommand{\eP}{\end{proposition}}
\newcommand{\bL}{\begin{lemma}}
\newcommand{\eL}{\end{lemma}}
\newcommand{\bT}{\begin{theorem}}
\newcommand{\eT}{\end{theorem}}
\newcommand{\bC}{\begin{corollary}}
\newcommand{\eC}{\end{corollary}} 
\newcommand{\eop}{\hfill \rule{2mm}{2mm}}
\newcommand{\pf}{\noindent{\bf Proof.~}}
\newcommand{\PD}{\text{proj} \dim}
\newcommand{\lra}{\longrightarrow}
\newcommand{\hra}{\hookrightarrow}
\newcommand{\llra}{\longleftrightarrow}
\newcommand{\Lra}{\Longrightarrow}
\newcommand{\Llra}{\Longleftrightarrow}
\newcommand{\bE}{\begin{enumerate}}
\newcommand{\eE}{\end{enumerate}}
\newcommand{\Sets}{\underline{{\mathrm Sets}}}
\newcommand{\Sch}{\underline{{\mathrm Sch}}}
\newcommand{\ForMe}{\noindent\TCP{{\bf Remarks To Be Removed:~}}}
\newcommand{\pic}{The proof is complete.}
\newcommand{\tcp}{This completes the proof.}


\def\m{\mathfrak {m}}
\def\CA{\mathcal {A}}
\def\CB{\mathcal {B}}
\def\CP{\mathcal {P}}
\def\CC{\mathcal {C}}
\def\CD{\mathcal {D}}
\def\CE{\mathcal {E}}
\def\CF{\mathcal {F}}
\def\CE{\mathcal {E}}
\def\CG{\mathcal {G}}
\def\CH{\mathcal {H}}
\def\CI{\mathcal {I}}
\def\CJ{\mathcal {J}}
\def\CK{\mathcal {K}}
\def\CL{\mathcal {L}}
\def\CM{\mathcal {M}}
\def\CN{\mathcal {N}}
\def\CO{\mathcal {O}}
\def\CP{\mathcal {P}}
\def\CQ{\mathcal {Q}}
\def\CR{\mathcal {R}}
\def\CS{\mathcal {S}}
\def\CT{\mathcal {T}}
\def\CU{\mathcal {U}}
\def\CV{\mathcal {V}}
\def\CW{\mathcal {W}}
\def\CX{\mathcal {X}}
\def\CY{\mathcal {Y}}
\def\CZ{\mathcal {Z}}

\newcommand{\smallcirc}[1]{\scalebox{#1}{$\circ$}}
\def\BA{\mathbb {A}}
\def\BB{\mathbb {B}}
\def\BC{\mathbb {C}}
\def\BD{\mathbb {D}}
\def\BE{\mathbb {E}}
\def\BF{\mathbb {F}}
\def\BG{\mathbb {G}}
\def\BH{\mathbb {H}}
\def\BI{\mathbb {I}}
\def\BJ{\mathbb {J}}
\def\BK{\mathbb {K}}
\def\BL{\mathbb {L}}
\def\BM{\mathbb {M}}
\def\BN{\mathbb {N}}
\def\BO{\mathbb {O}}
\def\BP{\mathbb {P}}
\def\BQ{\mathbb {Q}}
\def\BR{\mathbb {R}}
\def\BS{\mathbb {S}}
\def\BT{\mathbb {T}}
\def\BU{\mathbb {U}}
\def\BV{\mathbb {V}}
\def\BW{\mathbb {W}}
\def\BX{\mathbb {X}}
\def\BY{\mathbb {Y}}
\def\BZ{\mathbb {Z}}

\newcommand{\TCP}{\textcolor{purple}}
\newcommand{\TCM}{\textcolor{magenta}}
\newcommand{\TCR}{\textcolor{red}}
\newcommand{\TCB}{\textcolor{blue}}
\newcommand{\TCG}{\textcolor{green}}

\def\SA{\mathscr {A}}
\def\SB{\mathscr {B}}
\def\SC{\mathscr {C}}
\def\SD{\mathscr {D}}
\def\SE{\mathscr {E}}
\def\SF{\mathscr {F}}
\def\SG{\mathscr {G}}
\def\SH{\mathscr {H}}
\def\SI{\mathscr {I}}
\def\SJ{\mathscr {J}}
\def\SK{\mathscr {K}}
\def\SL{\mathscr {L}}
\def\SN{\mathscr {N}}
\def\SO{\mathscr {O}}
\def\SP{\mathscr {P}}
\def\SQ{\mathscr {Q}}
\def\SR{\mathscr {R}}
\def\SS{\mathscr {S}}
\def\ST{\mathscr {T}}
\def\SU{\mathscr {U}}
\def\SV{\mathscr {V}}
\def\SW{\mathscr {W}}
\def\SX{\mathscr {X}}
\def\SY{\mathscr {Y}}
\def\SZ{\mathscr {Z}}

\def\bfA{{\bf A}}
\def\bfB{{\bf B}} 
\def\bfC{{\bf C}} 
\def\bfD{{\bf D}} 
\def\bfE{{\bf E}} 
\def\bfF{{\bf F}} 
\def\bfG{{\bf G}} 
\def\bfH{{\bf H}} 
\def\bfI{{\bf I}} 
\def\bfJ{{\bf J}} 
\def\bfK{{\bf K}} 
\def\bfL{{\bf L}} 
\def\bfM{{\bf M}} 
\def\bfN{{\bf N}} 
\def\bfO{{\bf O}} 
\def\bfP{{\bf P}} 
\def\bfQ{{\bf Q}} 
\def\bfR{{\bf R}} 
\def\bfS{{\bf S}} 
\def\bfT{{\bf T}} 
\def\bfU{{\bf U}} 
\def\bfV{{\bf V}} 
\def\bfW{{\bf W}} 
\def\bfX{{\bf X}} 
\def\bfY{{\bf Y}} 
\def\bfZ{{\bf Z}} 

\def\a{\mathfrak {a}}
\def\b{\mathfrak {b}}
\def\c{\mathfrak {c}}

\def\d{\mathfrak {d}}
\def\e{\mathfrak {e}}
\def\f{\mathfrak {f}}
\def\g{\mathfrak {g}}
\def\i{\mathfrak {i}}
\def\j{\mathfrak {j}}
\def\k{\mathfrak {k}}
\def\l{\mathfrak {l}}
\def\m{\mathfrak {m}}
\def\n{\mathfrak {n}}
\def\p{\mathfrak {p}}
\def\q{\mathfrak {q}}
\def\r{\mathfrak {r}}
\def\s{\mathfrak {s}}
\def\t{\mathfrak {t}}
\def\u{\mathfrak {v}}
\def\w{\mathfrak {w}}

\def\BIGoplus#1#2{{\stackrel{#1}{#2}}}
\maketitle
\noindent{\bf Abstract:}
Let $X=\spec{A}$ denote a 
regular affine scheme, over a field $k$, with $1/2\in k$ and $\dim X=d$. Let $P$ denote a projective $A$-module of rank $n\geq 2$. 
Let $\pi_0\left({\mathcal LO}(P)\right)$ denote the (Nori) Homotopy Obstruction set \cite{M2, MM2}, and $\widetilde{CH}^n\left(X, \Lambda^nP\right)$ 
denote the Chow Witt group, of Barge and Morel \cite{BM}. In this article, we define a natural (set theoretic) map
$$
\Theta_P: \pi_0\left({\mathcal LO}(P)\right) \lra \widetilde{CH}^n\left(X, \Lambda^nP\right)
$$
{\bf The main results are published in the book \cite{M23}, and this article remains unpublished here in arxiv.}

\section{Introduction} 

{\bf This is an update for the readers. }
The main results in this article are included in my recently published Book \cite{M23}.
This is a landmark paper on Obstruction theory of projective modules. Since the results are already available in the book \cite{M23} it will not make sense to publish this paper again. 
As usual, the paper remains here in Arxiv, for reference. 
Only problem left in this theory is the Agreement question (\ref{AgeeeQb}). 
It has been established that the Nori Homotopy obstruction detects splitting \cite{MM2}, under suitable conditions. Question remains whether Chow-Witt obstruction of Barge-Morel would also do the same. 
I would leave this problem for the Chow-Witt cohort to resolve.

In this article we establish a natural (set theoretic)  map  from (Nori) Homotopy Obstruction set \cite{M2, MM2},
 to the Chow Witt group obstructions \cite{BM}, 
for projective modules
to split off a free direct summand. We avoid repeating the extensive background comments given in \cite{MM1, MM2},
on this  set of problems on obstructions for projective modules. 
  To facilitate further discussions, in this introduction, let $X=\spec{A}$ denote a 
regular affine scheme, over a field $k$, with $1/2\in k$ and $\dim X=d$. Let $P$ denote a projective $A$-module, with   $rank(P)=n$.
The episode began, when  (around 1989) through some verbal communications M. V. Nori  posed the following Homotopy question. 

%
\begin{question}
[Homotopy Question]\label{homoConj}
Suppose $X=\spec{A}$ is a smooth affine variety, with $\dim X=d$. Let $P$ be a projective $A$-module of rank $n$ and $f_0:P\sur I$ 
be  a surjective 
homomorphism, onto an ideal $I$ of $A$. Assume $Y=V(I)$ is smooth with $\dim Y=d-n$.
Also suppose $Z=V(J)\subseteq \spec{A[T]}=X\times \BA^1$ is a smooth subscheme, such that $Z$ intersects $X\times 0$ transversally in $Y\times 0$.
Now, suppose that $\varphi: P[T] \sur \frac{J}{J^2}$ is a surjective map such that $\varphi_{|T=0}=f_0\otimes \frac{A}{I}$.  
Then the question is, whether there is a
surjective 
map $F:P[T]\sur J$ such that (i) $F_{|T=0}=f_0$ and (ii) $F_{|Z}=\varphi$.
{\rm Assume $2n\geq d+3$.}
\end{question}
 Refer to \cite{MM2, Mu, M2} and others for slightly varying versions of this question (\ref{homoConj})).  
 The Homotopy question (\ref{homoConj}) was settled  affirmatively  
by Bharwadekar-Keshari \cite{BK}, when $A$ is essentially smooth over an infinite perfect field $k$, and $2n\geq d+3$. The proof used the
corresponding result  
in the local case in \cite{MV}, along with results in \cite{M2}.

\noindent{\bf Euler class groups:}  
Subsequent to 
 (\ref{homoConj}), Nori also communicated (1992) a definition of  Euler class groups to house obstructions for splitting,
in the case $d=n$ (see \cite{MR}). 
Given integers $0\leq n\leq d$ and invertible modules ${\CL}$, Nori's definition was expanded (\cite{BS1, BS2, MY}) to define 
 Euler class groups $E^n(A, {\CL})$. These were defined ideal theoretically, using generators and relations.
However, only when $n=d$, the Euler class groups worked well enough.  
 For projective $A$-modules $P$ of 
 rank $d$, an Euler class  $e(P) \in E^d(A, \Lambda^d P)$ was defined. It was conjectured and proved \cite{BS2, BS3} that 
\begin{equation}\label{272AsinEgrgr}
 e(P)=0 \Llra P \cong Q\oplus A.
\end{equation} 
 For projective $A$-modules $P$, with $rank(P)=n< d$, attempts to define Euler classes $e(P)$ failed. 
%


This set of ideas of Nori ({\it along with other ongoing activities,} e.g. \cite{BS1, BS2, BS3, Mu, M2}) received an added significance in 2000, with the introduction of 
 Chow Witt groups, by Barge and Morel \cite{BM, Mo, F1}, to house such a possible obstruction. 
 ({\it Chow Witt groups are also known as Oriented Chow grous.})
%
For  integers $0\leq n\leq d$, and  line bundles ${\CL}$ on $X$,   groups $\widetilde{CH}^n\left(X, {\CL}\right)$, to be called the 
Chow Witt groups,
 were introduced. Further, 
 with $rank(P)=n$,
  an obstruction class $\varepsilon_{CW}(P) \in  \widetilde{CH}^n\left(X, \Lambda^n P\right)$, to be called the Chow Witt obstruction,
was defined. As in the case of Euler class groups (\ref{272AsinEgrgr}), when $n=d$, it was conjectured \cite{BM}
 and proved that \cite[Theorem 8.14]{Mo}
\begin{equation}\label{323MorelEpCW}
\varepsilon_{CW}(P) =0  \Llra P \cong Q\oplus A
\end{equation} 
when $A$ is smooth over a perfect field $k$.
 By then, it started appearing not so promising, that Euler class groups $E^n(A, \Lambda^n P)$ may  be able to house such obstructions, for splitting,
 when $n=rank(P)\leq d-1$. It remains an open question,  whether Chow Witt obstruction $\varepsilon_{CW}(P) \in  \widetilde{CH}^n\left(X, \Lambda^n P\right)$,
 would detect splitting, under certain reasonable conditions, when $n=rank(P)\leq d-1$.

 While it was not explicitly articulated, it is obvious  that there is an ingrained homotopy relation  in the Homotopy question (\ref{homoConj}).
For this reason, a local $P$-orientation is defined to be a pair $(I, \omega)$, where $I\subseteq A$ is an ideal 
and $\omega: P \sur \frac{I}{I^2}$ is a surjective map. Let ${\mathcal LO}(P)$ denote the set of all local $P$-orientations.
By substituting $T=0, 1$, one obtains two maps
 $$
 \diagram
  {\mathcal LO}(P) & {\mathcal LO}(P[T])\ar[r]^{T=1} \ar[l]_{T=0} & {\mathcal LO}(P) \\
  \enddiagram
 $$
This induces chain equivalence relations $\sim$ on ${\mathcal LO}(P)$. 
The homotopy obstruction set $\pi_0\left( {\mathcal LO}(P)\right)$
was defined to be the set of all equivalence classes.
Further, an obstruction class $\varepsilon_H(P)\in \pi_0\left({\mathcal LO}(P)\right)$ is defined. It was also established \cite{MM2}
that 
\begin{equation}\label{MMDuiResult}
\varepsilon_H(P)\quad {\rm is~neutral} \quad \Llra \quad P \cong Q\oplus A
\end{equation}
 when $2n\geq d+3$ and 
$A$ is essentially smooth over a perfect field. This is clearly an improvement over the results (\ref{272AsinEgrgr}, \ref{323MorelEpCW}), which works only when the rank $n=d$.
The possibility to use the homotopy relations ingrained in  (\ref{homoConj}) to construct a house $\pi_0\left({\mathcal LO}(P)\right)$ for obstructions $\varepsilon_H(P)$ 
was considered
only recently \cite{MM2}. 
However, $\pi_0\left({\mathcal LO}(P)\right)$ is an invariant of $P$ itself, 
while Chow Witt groups $\widetilde{CH}^n\left(X, {\CL}\right)$ are  invariants of $X$ (and of the determinant). 
Further, as invariants, Chow Witt groups $\widetilde{CH}^n\left(X, {\CL}\right)$ are very similar to Chow groups $CH^n(X)$ \cite{Fu}. 
%
We summarize the some of the above:
\bE
\item  The homotopy obstruction $\varepsilon_H(P) \in \pi_0\left({\mathcal LO}(P)\right)$ is well defined and it detects splitting (\ref{MMDuiResult}), 
\item  The Chow Witt obstruction $\varepsilon_{CW}(P) \in \widetilde{CH}^n\left(X, \Lambda^n P\right)$ is well defined. Only when $n=d$, it is known to detect 
splitting (\ref{323MorelEpCW}). 
\item  The Euler class $e(P)\in E^d(A, \Lambda^dP)$ is defined only when $n=d$ and it detects splitting in this case. For $n\leq d-1$, it 
is unlikely that an obstruction (Euler) class $e(P)$ in   $E^n(A, \Lambda^nP)$ would be definable. Further, 
in \cite{MM2}, for any projective module $P$, an Euler class group $E(P)$ is defied. 
This group is generated by the set ${\mathcal LO}^n(P)$ of all $P$ orientations 
$(I, \omega)$ with $height(I)=n$; modulo the global orientations. Further, $E^n(A, {\CL})=E({\CL}\oplus A^{n-1})$. 
The Euler class groups $E(P)$ works very well when $n=d$ \cite{MM1}. While the definition of $E(P)$ is very natural, it fails decisively when $n\leq d-1$.
Under some stringent conditions, there is a natural map $E(P) \sur \pi_0\left({\mathcal LO}(P)\right)$ \cite{MM2}. 
\eE 


At this stage it is fairly transparent that Homotopy obstructions of Nori \cite{M2, MM2},
 and the Chow Witt obstructions of Barge-Morel \cite{BM} must converge, some way.
This is precisely what we
respond to, in this article, by establishing a natural map (set theoretic) 
\begin{equation}\label{389tHetaP}
\Theta_P: \pi_0\left({\mathcal LO}(P) \right) \lra \widetilde{CH}^n\left(X, \Lambda^n P\right)\quad {\rm with}\quad 
\Theta_P\left(\varepsilon_{H}(P)\right) = \varepsilon_{CW}(P^*).
\end{equation}
We establish this for all projective modules $P$, with $2\leq n=rank(P) \leq d$. While it was established in \cite{MM2} that $\pi_0\left({\mathcal LO}\right)(P)$
has a additive structure, when $2n\geq d+2$, there is no such well defined structure outside this range of $n$. Therefore, the map 
$\Theta_P$, defined above, would be 
a set theoretic map only, in general, and respects additivity when $2n\geq d+2$. 
Further, recall that the obstruction class $\varepsilon_{H}(P) \in \pi_0\left({\mathcal LO}(P)\right)$ is defined,
for all $n:=rank(P)$. The obstruction class $\varepsilon_{H}(P)$ detects splitting properties of $P$, under the conditions stated above (\ref{MMDuiResult}).
%
With this result in mind,  the following natural question emerges. 
\begin{question}
[Agreement Question]\label{AgeeeQb} 
{\rm 
Whether the map $\Theta_P$ is injective? 
}
\end{question}
If and when the answer to (\ref{AgeeeQb}) is affirmative,
Chow Witt obstructions $\varepsilon_{CW}(P)$ would detect splitting, under the same hypotheses above (\ref{MMDuiResult}).

This  set of problems was referred to  as the "Homotopy Program"  by this author (e.g. \cite{MS}).
 More precise outline of the program was given in \cite{MM2}. The possibility of the existence of
  a map $\Theta_P$, as above, was mentioned as a part of the program \cite[{\bf Part 2}, pp. 173]{MM2}, 
 which we accomplish fully 
 in this article. 
 
 We briefly comment on the place of ${\BA}^1$-homotopy in this program \cite[{\bf Part 2}, pp. 173]{MM2}. 
  We denote the ${\BA}^1$-homotopy category by ${\CH}(k)$. 
  Perhaps, this stems from variety of spheres $S^{p, q}$ considered in ${\CH}(k)$ \cite[pp.111]{MVv}, and 
  $Q_{2n}=\spec{\frac{k[X_1, \ldots, X_n; Y_1, \ldots, Y_n, Z]}{\left(\sum_{i=1}^nX_iY_i+Z(Z-1)\right)}}\cong ({\BP}^1_k)^{\wedge n}\cong S^{2n, n}$  being one of them
  \cite[Rem. 6.42]{Mo}, \cite[Thm. 2.2.5]{ADF}. 
  Further, $Q_{2n}$ is ${\BA}^1$-naive \cite[Thm. 1.1.1]{AF}, in the sense 
\begin{equation}\label{432UnotoIsos}
  \pi_{{\BA}^1}(X, Q_{2n}) \iso Mor_{{\CH}(k)}\left(X, Q_{2n}\right)\quad \forall~{\rm smooth~affine~schemes}~ X=\spec{A},
\end{equation}
where, following the Suslin-Voevodsky construction, $\pi_{{\BA}^1}(X, Q_{2n})$ denotes the usual set defined by homotopy (see \cite[pp. 199]{Mo}).
Other notations were used in \cite{MM1, MM2, MM3}, for $\pi_{{\BA}^1}(X, Q_{2n})$.
Combining this with results in \cite[Lem. 2.9, Thm. 7.3]{MM2}, when $2n\geq d+2$,
we have maps and isomorphisms
\begin{equation}\label{431Isos}
E^n(A, A)= E(A^n) \lra \pi_0({\mathcal LO}(A^n) \cong   \pi_{{\BA}^1}(X, Q_{2n}) \iso Mor_{{\CH}(k)}\left(X, Q_{2n}\right) 
\end{equation}
The second arrow is also an isomorphism, when $A$ is essentially smooth, $k$ is an infinite perfect field, and $2n\geq d+3$ \cite[Thm. 1.4]{MM1}.
 Combining with  (\ref{389tHetaP}), we obtain a natural map
\begin{equation}\label{440Did}
 Mor_{{\CH}(k)}\left(X, Q_{2n}\right) \lra \widetilde{CH}^n\left(X, A\right) \quad {\rm if}~2n\geq d+3.
\end{equation}
This map (\ref{440Did}) was also obtained in \cite[Thm. 1]{AF}, when $2n\geq d+2$, which is a particular case of (\ref{389tHetaP}), with $P=A^n$. 
%
%
Given these isomorphisms (\ref{431Isos}), the ${\BA}^1$-homotopy invariant $Mor_{{\CH}(k)}\left(X, Q_{2n}\right)$ does not seem to be of any additional help, from our perspective of 
obstruction theory. However, as indicated in \cite[Appendix. A]{MM2}, the obstructions $\pi_0({\mathcal LO}(P))$ should have a ${\BA}^1$-homotopy interpretation, 
which would be of  interest, by its own rights. Such an interpretation may also throw some light on the question (\ref{AgeeeQb}). 

As is \cite{BM}, our work is fully reliant on the methods of classical (commutative) algebra, and is independent of the methods of ${\BA}^1$-homotopy.

%

We comment on the organization of the article. In section \ref{recapSec}, we provide some background on Homotopy obstruction, mainly from \cite{MM2, MM1}.
In section \ref{simformSec}, we associate a symmetric form $\Phi(I, \omega)$ to certain  representatives  $(I, \omega)$  of the elements 
$\left[(I, \omega)\right] \in \pi_0\left({\mathcal LO}(P) \right)$. In section \ref{BckCWSec}, we provide some preliminaries on Chow Witt groups. In section \ref{ConvSec},
we establish the map $\Theta_P$.

\noindent{\bf Acknowledgement.} {\it This author is thankful to Marco Schlichting for many valuable discussions.}

\section{Preliminaries on Homotopy Obstructions} \label{recapSec}  
Throughout this article $A$ will denote a noetherian commutative ring, with $\dim A=d$, and
$A[T]$ will denote the polynomial
ring in one variable $T$. We  assume $1/2\in A$ and $\spec{A}$ is connected. 
For an $A$-module $M$, denote
$M[T]:=M\otimes A[T]$. Likewise, for a homomorphism $f:M\lra N$ of  $A$-modules, $f[T]:=f\otimes A[T]$. 
Our main results would assume $A$ is a regular ring, containing a field $k$, with $1/2\in k$.
Further, $P$ will denote a projective $A$-module with $rank(P)=n$, and $2 \leq n \leq d$. Denote $\left|P\right|:=\Lambda^nP$, the determinant of $P$. 

%
%
For a such a projective $A$-module $P$, as above, the set $\pi_0\left({\mathcal LO}(P) \right)$ of equivalence classes of homotopy obstructions 
was defined in \cite{MM2}. We 
recall some of the  essential elements of the definition of, and alternative descriptions of $\pi_0\left({\mathcal LO}(P)\right)$ from \cite{MM1}.

\bD\label{defpizero2Oct}{\rm 
Let $A$ be a noetherian commutative  ring, with $\dim A=d$ and $P$ be a projective $A$-module with $rank(P)=n$.
 By a {\bf local $P$-orientation}, we mean 
a pair $(I, \omega)$ where $I$ is an ideal of $A$ and $\omega:P \sur \frac{I}{I^2}$ is a surjective homomorphism.
We will use the same notation $\omega$ for the map $\frac{P}{IP} \sur \frac{I}{I^2}$, induced by $\omega$. 
 A local {\bf local $P$-orientation} will simply be 
 referred to as a {\bf local orientation}, when $P$ is understood.
Denote 
\begin{equation}\label{4LOetc}
\left\{
\begin{array}{l}
{\mathcal LO}(P)=\left\{(I, \omega): (I, \omega)~{\rm is~a~local}~P~{\rm orientation} \right\}\\
\widetilde{{\mathcal LO}}(P)=\left\{(I, \omega)\in {\mathcal LO}(P): height(I)\geq n  \right\}\\
{\mathcal LO}^{n}(P)=\left\{(I, \omega)\in {\mathcal LO}(P): height(I)= n  \right\}\\
{\rm Note}  ~\widetilde{{\mathcal LO}}(P)= {\mathcal LO}^{n}(P) \cup\{(A, 0)\}\\
\end{array}
\right.
\end{equation}
%
%

For $(I_0, \omega_0), (I_1, \omega_1)\in {\mathcal LO}(P)$, 
we write 
$(I_0, \omega_0)\sim_H (I_1, \omega_1)$, if there is an $(I, \omega)\in {\mathcal LO}(P[T])$ such that  $(I(0), \omega(0))=(I_0, \omega_0)$ and 
$(I(1), \omega(1))=(I_1, \omega_1)$. 
 In this case, we say $(I_0, \omega_0)$ is homotopic to $(I_1, \omega_1)$.
The realtion $\sim_H$ generates a chain equivalence relation $\sim$ on ${\mathcal LO}(P)$, which we call the chain homotopy relation.
The (Nori) {\bf Homotopy obstruction} set $\pi_0\left({\mathcal LO}(P) \right):=\frac{{\mathcal LO}(P)}{\sim}$ is
defined to be the set of equivalence classes of elements in ${\mathcal LO}(P)$.
%
  So, we have a push forward diagram:
\begin{equation}\label{pizeroLOpush}
\diagram
{\mathcal LO}(P[T])\ar[r]^{T=0}\ar[d]_{T=1} &{\mathcal LO}(P)\ar[d]\\
{\mathcal LO}(P)\ar[r] & \pi_0\left({\mathcal LO}(P)\right)\\ 
\enddiagram 
\qquad {\rm in}\quad \underline{Sets}.
\end{equation}

Similarly, the homotopy relation $\sim_H$ on $\widetilde{{\mathcal LO}}(P)$ leads to a chain equivalence relation $\sim$ on $\widetilde{{\mathcal LO}}(P)$.  
Note, due to moving lemma arguments (basic element theory), it would not make any difference,
whether $\sim_H$ 
 is defined by using  $P[T]$-orientations $(I, \omega)$ in $\widetilde{{\mathcal LO}}(P[T])$ or in ${\mathcal LO}(P[T])$. 
Define  
$$
\pi_0\left(\widetilde{{\mathcal LO}}(P) \right):=\frac{\widetilde{{\mathcal LO}}(P)}{\sim}.
$$
 Note, a push forward diagram, as in in (\ref{pizeroLOpush}),
to define $\pi_0\left(\widetilde{{\mathcal LO}}(P) \right)$, would not be possible, because the substitution $T=0, T=1$ do not behave too well, in this case. 
%
%
%
%
%
However,  there is natural map
$$
\varphi: \pi_0\left(\widetilde{{\mathcal LO}}(P)\right) \lra \pi_0\left({\mathcal LO}(P)\right)
$$
}
\eD
%
The following is from \cite{MM2}.
\bP\label{transIso}{\rm 
Let $A$ and $P$ be as in (\ref{defpizero2Oct}). Then, the map $\varphi$ is surjective.
Assume further that $A$ is a regular ring  containing a field $k$, with $1/2\in k$.
Then $\sim$ is  an equivalence relation on $\widetilde{{\mathcal LO}}(P)$. 
Moreover, 
 $\varphi$ is a bijection.
}%
\eP
%

%
\begin{remark}\label{moviongRem}{\rm 
With notations as in (\ref{defpizero2Oct}) the following are some useful observations.
\bE
\item Suppose $(I_0, \omega_0), (I_1, \omega_1)\in \widetilde{{\mathcal LO}}(P)$ and $(I_0, \omega_0)\sim (I_1, \omega_1)$. By defintion there is homotopy
$H(T)=(I, \omega)\in {\mathcal LO}(P[T])$ such that $H(0)= (I_0, \omega_0)$ and $H(1)=(I_1, \omega_1)$.  
By moving Lemma arguments, 
similar to \cite[Lemma 4.5]{MM2}, we can assume that $H(T)\in \widetilde{{\mathcal LO}}(P[T])$.
%
\item If $A$ is  Cohen Macaulay, then $\widetilde{{\mathcal LO}}(P)$ is in bijection with the set
$$
\left\{(I, \omega) \in {\mathcal LO}^n(P): \omega: \frac{P}{IP} \iso \frac{I}{I^2}~{\rm is~an~isomorphism}\right\}
\cup \left\{(A, 0)\right\}
$$
Moreover, for $(I, \omega) \in {\mathcal LO}^n(P)$, $I$ is a local complete intersection ideal. 
\eE
}
\end{remark}


\vspace{3mm}
We recall the definition of the neutral element and the obstruction class.
\bD\label{neutrlObsCls}{\rm 
Use the  notations as above (\ref{defpizero2Oct}). The neutral element  ${\bf e}_1\in \pi_0\left({\mathcal LO}(P)\right)$ is defined to be the image of $(A, 0)$.
The Nori Homotopy obstruction $\varepsilon_H(P) \in \pi_0\left({\mathcal LO}(P)\right)$ is defined to be the image of $(0, 0)$. 

Now assume $A$ is a Cohen Macaulay ring and  $(I, \omega) \in \widetilde{{\mathcal LO}}(P)$. Suppose $\omega$ lifts to a surjective map 
$f: P \sur I$. It follows 
$$
\varepsilon_H(P)= image (I, \omega) \in \pi_0\left({\mathcal LO}(P)\right)
$$
}
\eD


\section{The symmetric form}\label{simformSec}
The essence of the arguments in the this section, can be traced back  to the following theorem of Altman and Kleiman \cite[Theorem 4.5]{AK}.

\bT\label{keyAltmanKleiman}{\rm 
Suppose $A$ is a commutative noetherian ring and $I$ is a locally complete intersection ideal, with $height(I)=n$. 
Suppose ${\CL}$ is an invertible $A$-module. 
Then, there is a natural isomorphism
$$
\chi: Hom\left(\Lambda^n \frac{I}{I^2}, \frac{{\CL}}{I{\CL}} \right)  \iso Ext^n\left(\frac{A}{I}, {\CL} \right) 
$$
For convenience, recall other notations:
$$
Hom\left(\Lambda^n \frac{I}{I^2}, \frac{{\CL}}{I{\CL}} \right)=:{\CL}\left|\frac{I}{I^2}\right|^{-1} 
$$
}
\eT

\vspace{3mm}
For our purpose, the following formulations would be helpful. 
\bL\label{79Alklein}{\rm 
Let $A$ be a noetherian commutative ring and $P$ be a projective $A$-module with $rank(P)=n$. 
Let $I$ be a locally complete 
intersection ideal with $height(I)=n$. Suppose $\omega:\frac{P}{IP} \iso \frac{I}{I^2}$ is an isomorphism, and
let $\left|\omega\right|:\left|\frac{P}{IP}\right|  \iso \left|\frac{I}{I^2}\right|$ denote the determinant of $\omega$,
where $\left|\frac{P}{IP}\right|:=\Lambda^n\frac{P}{IP}$,  $\left|\frac{I}{I^2}\right|:=\Lambda^n\frac{I}{I^2}$ denote the determinants.

Let $f:P\sur I$ be a surjective lift of $\omega$, and consider its Koszul complex. 
Then, for any finitely generated $A$-module $M$, the sequence
\begin{equation}\label{592AltmanSeq}
\diagram 
Hom\left(\Lambda^{n-1}P, M\right) \ar[r]^{d^*_n} & Hom\left(\Lambda^{n}P, M\right)\ar[r]^{\varphi} & Hom\left(\Lambda^n\frac{P}{IP}, \frac{M}{IM} \right) \ar[r] & 0\\
\enddiagram 
\end{equation} 
is exact. Consequently, the following maps are isomorphisms:
\begin{equation}\label{604NatItta}
\diagram 
Hom\left(\left|\frac{P}{IP}\right|, \frac{M}{IM} \right)  \ar[r]^{\quad\iota(f)}_{\quad\sim} &Ext^n\left(\frac{A}{I}, M\right) \\ 
Hom\left(\left|\frac{I}{I^2}\right|, \frac{M}{IM} \right)\ar[u]^{\widetilde{\left|\omega\right|}}_{\wr} \ar@/_/@{-->}[ru]_{\chi(\omega)} &\\
\enddiagram
 {\rm where}
\left\{\begin{array}{l}
\iota(f)~{\rm is ~the ~  identification}\\
\widetilde{\left|\omega\right|} \quad{\rm sends}~\lambda\mapsto \lambda \left|\omega\right|\\ 
\chi(\omega) =\iota(f)\widetilde{\left|\omega\right|} \\
\end{array}
\right.
\end{equation}
Further, the map $\chi(\omega)$ depends only on $\left|\omega\right|$, and is independent of the lift $f$ of $\omega$.

}
\eL 
\pf Note there is a natural surjective map $\varphi$, by reduction modulo $I$. Further, the composition is zero. We check the exactness of (\ref{592AltmanSeq}) locally.
So, we can assume $P=\oplus_{i=1}^nAe_i$ is free, and let $f(e_i)=x_i$. Suppose $\varphi(\lambda)=0$ for fome $\lambda:Ae_1\wedge \cdots \wedge e_n \lra M$. 
So, 
$$
\lambda(e_1\wedge \cdots \wedge e_n)=\sum_{i=1}^n x_im_i \qquad {\rm for~some}\quad m_i\in M
$$
Define
$$
\widetilde{\lambda}: \Lambda^{n-1}P \lra M \quad {\rm by}\quad \widetilde{\lambda}(e_1\wedge \cdots \wedge \hat{e}_i \wedge\cdots \wedge e_n)=
(-1)^{i-1}m_i
$$
Then, $d_n^*(\widetilde{\lambda})=\lambda$. This establishes that that sequence (\ref{592AltmanSeq}) is exact.

The map $\varphi$ is independent of the lift $f$. So, by definition $\chi(\omega)$ depends only on $\left|\omega\right|$.
 \pic $\eop$

\vspace{3mm}
\begin{remark}\label{twoExt644}{\rm 
Recall, as always, $Ext^n\left(\frac{A}{I}, M\right)$ in (\ref{79Alklein}) is defined only up to isomorphism and  
 the natural identification $\iota$ in (\ref{604NatItta}) needs to be understood with care. Suppose $A$, $I$, $P$, $M$ be as in (\ref{79Alklein}), and let 
 $f, g:P \sur I$ be two surjective maps, ({\it without any reference to} $\omega$). Then, corresponding to Koszul complexes of $f$, $g$ are related as follows
 $$
 \diagram
0\ar[r] & \Lambda^nP \ar[r]^{d_n(f)}\ar[d]_{\det(\Delta)}   & \Lambda^{n-1}P   \ar[r] \ar[d]& \cdots \ar[r] & P\ar[r]^f\ar[d]^{\Delta} & A \ar[r]\ar[d]^1 & \frac{A}{I} \ar[r] \ar[d]^1& 0\\
0\ar[r] & \Lambda^nP \ar[r]_{d_n(g)}  \ar[r] & \Lambda^{n-1}P \ar[r] \ar[r] & \cdots \ar[r] & P\ar[r]_g & A \ar[r] & \frac{A}{I} \ar[r] & 0\\
 \enddiagram 
 $$ 
 By (\ref{592AltmanSeq}), we have two representations $Hom\left(\Lambda^n\frac{P}{IP}, \frac{M}{IM} \right)\iso Ext^n\left(\frac{A}{I}, M\right)$, 
 related, as follows
$$
\diagram 
Hom\left(\Lambda^{n-1}P, M\right) \ar[r]^{d(f)^*_n}\ar[d] & Hom\left(\Lambda^{n}P, M\right)\ar[r]^{\varphi}\ar[d]^{\det(\Delta)} 
& Hom\left(\Lambda^n\frac{P}{IP}, \frac{M}{IM} \right) \ar[r]\ar[d]^{\overline{\det(\Delta)}}  & 0\\
Hom\left(\Lambda^{n-1}P, M\right) \ar[r]_{d(g)^*_n} & Hom\left(\Lambda^{n}P, M\right)\ar[r]_{\varphi} & Hom\left(\Lambda^n\frac{P}{IP}, \frac{M}{IM} \right) \ar[r] & 0\\
\enddiagram 
$$
where $\overline{\det(\Delta)}=\Delta ~{\rm mod}~I$ is an unit. One can summarize, that the diagram 
$$
\diagram 
 Hom\left(\left|\frac{P}{IP}\right|, \frac{M}{IM} \right)\ar[d]_{\overline{\det(\Delta)}}^{\wr} \ar[r]^{\iota(f)} & Ext^n\left(\frac{A}{I}, M\right)\\
Hom\left(\left|\frac{P}{IP}\right|, \frac{M}{IM} \right)  \ar@/_/[ur]_{\iota(g)}&\\
\enddiagram 
\qquad \qquad {\rm commutes.}
$$
}
\end{remark} 

\vspace{3mm}
The following is a further relaxed version of (\ref{79Alklein}), which suits our purpose. 
\bC\label{79AlTNori}{\rm 
Let $A$ be a noetherian commutative ring and $P$ be a projective $A$-module with $rank(P)=n$. 
Let $I$ be a locally complete 
intersection ideal with $height(I)=n$, and $\omega:\frac{P}{IP} \iso \frac{I}{I^2}$ be an isomorphism. 
Let $f:P\sur IJ$ be a surjective lift of $\omega$, where $J$ is an ideal with $I+J=A$.

%
Then, for any finitely generates $A$-module $M$, 
the following maps are isomorphisms:
\begin{equation}\label{698Exi6}
\diagram 
Hom\left(\left|\frac{P}{IP}\right|, \frac{M}{IM} \right)  \ar[r]^{\quad\iota(f)}_{\quad\sim} &Ext^n\left(\frac{A}{I}, M\right) \\ 
Hom\left(\left|\frac{I}{I^2}\right|, \frac{M}{IM} \right)\ar[u]^{\widetilde{\left|\omega\right|}}_{\wr} \ar@/_/@{-->}[ru]_{\chi(\omega)} &\\
\enddiagram
 {\rm where}
\left\{\begin{array}{l}
\iota(f)~{\rm is ~the ~  identification}\\
\widetilde{\left|\omega\right|} \quad{\rm sends}~\lambda\mapsto \lambda \left|\omega\right|\\ 
\chi(\omega) =\iota(f)\widetilde{\left|\omega\right|}\\
\end{array}\right.
\end{equation}
Further, the map $\chi(\omega)$ depends only on $\left|\omega\right|$, and is independent of the lift $f$ of $\omega$, or the ideal $J$.

}
\eC  
\pf Let $s+t=A$ with $s\in I$ and $t\in J$. Then, the natural map $Ext^n\left(\frac{A}{I}, M\right)\iso Ext^n\left(\frac{A}{I}, M_t\right)$ is isomorphism.
Further, $f_t:P_t \sur I_t$ is a surjective map. It follows from (\ref{79Alklein}), 
the sequence
\begin{equation}\label{83WithJay}
\diagram 
Hom\left(\Lambda^{n-1}P_t, M_t\right) \ar[r]^{d^*_n} & Hom\left(\Lambda^{n}P_t, M_t\right)\ar[r]^{\varphi} & Hom\left(\Lambda^n\frac{P}{IP}, \frac{M}{IM} \right) \ar[r] & 0\\
\enddiagram 
\end{equation} 
is exact. The rest follows as in (\ref{79Alklein}).
\pic $\eop$



\vspace{3mm}

In the rest of this section, we further elaborate  all these (\ref{698Exi6}), to associate symmetric isomorphism $(I, \omega) \mapsto \Phi(I, \omega)$.

%
%
Suppose $X=\spec{A}$ is Cohen Macaulay scheme, 
with $\dim X=d$. Recall, for integers $0\leq n \leq d$, 
in \cite{M1}, the subcategory 
$C{\BM}^n(X)\subseteq Coh(X)$, was defined to be the full subcategory of objects
\begin{equation}\label{587CBMnX}
C{\BM}^n(X)=\left\{{\CF}\in Coh(X): grade({\CF})={\PD}({\CF})=n \right\}
\end{equation} 
where ${\PD}({\CF})$ denotes the locally free dimension of ${\CF}$. Given an invertible sheaf ${\CL}$, on $X$, the association 
${\CF} \mapsto {\CF}^{\vee}:={\SE}xt^n\left({\CF}, {\CL}\right)$ is a duality in $C{\BM}^n(X)$. This endows $C{\BM}^n(X)$ with a structure 
of an exact category with duality, to be denoted by  $C{\BM}^n\left(X, {\CL}\right)$. The following is a key definition for our purpose.

\bD\label{74DefPhira}{\rm
Now assume  $A$ is Cohen Macaulay ring, with $\dim A=d$. Let $P$ be a projective $A$-module with $rank(P)=n\leq d$. 
With the notations as in  section \ref{recapSec}, let $(I, \omega) \in {\mathcal LO}^n(P)$.  As before, $\left|\omega\right|: \left|\frac{P}{IP}\right| \iso
\left|\frac{I}{I^2}\right|$
denote the determinant map. With $M=\left|P \right|$, the diagram  (\ref{698Exi6}) reduces to:
\begin{equation}\label{779TryShivPhi}
\diagram 
End\left(\left|\frac{P}{IP}\right|\right) 
\ar[r]^{~\iota(f)}_{~\sim} &Ext^n\left(\frac{A}{I}, \left|P \right|\right) \\ 
 Hom\left(\left|\frac{I}{I^2}\right|, \left|\frac{P}{IP} \right| \right)\ar[u]^{\widetilde{\left|\omega\right|}}_{\wr} \ar@/_/@{-->}[ru]_{\chi(\omega)}&\\
\enddiagram
\end{equation}
Therefore,
$$
\chi(\omega)(\left|\omega\right|^{-1})=\iota(f)(1_{\left|\frac{P}{IP}\right|}) \in Ext^n\left(\frac{A}{I}, \left|P\right|\right) 
\quad {\rm  is ~a~free~generator.} 
$$
Define the  isomorphism 
\begin{equation}\label{AltKleiExt}
\f_I: \frac{A}{I} \iso Ext^n\left(\frac{A}{I}, \left|P\right|\right)~~~ {\rm by} \quad \f_I(\overline{a}) := a\iota(f)(1_{ \left|\frac{P}{IP} \right| })=\chi(\omega)(\left|\omega\right|^{-1})~~ \forall a\in A
\end{equation}

\bE
\item In other word, with identification $Ext^n\left(\frac{A}{I}, \left|P\right| \right) = End\left(\left|\frac{P}{IP}\right|\right)$, 
 the isomorphism $\f_I:\frac{A}{I} \iso End\left(\left|\frac{P}{IP}\right|\right)$,
sends $1\mapsto 1_{\left|\frac{P}{IP}\right|}$. 
\item Since  $Ext^r\left(\frac{A}{I}, \left|P\right| \right)=0,~\forall r\neq n$, it follows $\f_I$ is a symmetric isomorphism in $C{\BM}^n\left(X, \left|P\right|\right)$.

\item\label{810ItemPhiIomeg} We let $\Phi(I, \omega):=\left(\frac{A}{I}, \f_I \right)$ denote the symmetric space in $C{\BM}^n(X,\left|P\right|)$. 
  Further, for $(A, 0) \in \widetilde{{\mathcal LO}}(P)$, let $\Phi(A, 0):={\bf 0}$ be the trivial symmetric space in $C{\BM}^n(X, \left|P\right|)$. 
\item Let ${\bf M}\left(C{\BM}^n\left(X\right)\right):={\bf M}\left(C{\BM}^n\left(X, \left|P\right|\right)\right)$ denote the monoid of the isometry classes of 
symmetric spaces in $C{\BM}^n\left(X, \left|P\right|\right)$. We will continue to use the same notation $\Phi(I, \omega)$ for the symmetric space in 
$C{\BM}^n\left(X, \left|P\right|\right)$, and its isometry class in ${\bf M}\left(C{\BM}^n\left(X, \left|P\right|\right)\right)$. 
Thus, the association $(I, \omega) \mapsto \Phi(I, \omega)$, defines a set theoretic map:
\begin{equation}\label{801MoAss}
\Phi: \widetilde{{\mathcal LO}}(P) \lra {\bf M}\left(C{\BM}^n\left(X, \left|P\right|\right)\right)
\end{equation}

\eE 

}
\eD 
For future correspondences with the derived categories, we further analyze this symmetric space $\Phi(I, \omega)$.

 
\begin{remark}\label{DefPhiOmega231}{\rm 
Let $X=\spec{A}$ be an affine Cohen Macaulay scheme with $\dim X=d$.  
 Let $P$ be a projective $A$-module with $rank(P)=n$. 
 For any $A$-module $M$, write $M^*=Hom(A, \left|P\right|)$.
 Throughout, on the category ${\SP}(X)$ of finitely generated projective $A$-modules, we work with the duality $Q \mapsto Q^*$. Further, 
 duality we consider in other associated categories (e.g. $C{\BM}^n(X)$), would be induced by this duality on ${\SP}(X)$.
 
%
%
Let $(I, \omega)\in {\mathcal LO}^n(P)$.  
We have $\frac{A}{I} \in C{\BM}^n(X)$ and
 $\omega:\frac{P}{IP}\iso \frac{I}{I^2}$ is an isomorphism.  Consider a surjective lift $f$ of $\omega$, as in the diagram
 $$
 \diagram
 P \ar@{->>}[r]^{f} \ar@{->>}[d]&I \cap J\ar@{->>}[d]\\
 \frac{P}{IP}\ar[r]_{\omega}^{\sim}& \frac{I}{I^2}
 \enddiagram \quad {\rm with}~I+J=A,~height(J)\geq n
 $$
 Let $\omega_J:\frac{P}{JP} \iso \frac{J}{J^2}$, and $\omega_{IJ}:\frac{P}{IJP} \iso \frac{IJ}{I^2J^2}$ denote the isomorphism induced by $f$. Note 
 $\omega_{IJ}=\omega \oplus \omega_J$. 
   There is a natural map from the   Koszul complex of $f$ to its dual, as follows:
\begin{equation}\label{KoszulIsoMetry250}
 \diagram 
0\ar[r] & \Lambda^nP \ar[d]_{\varphi_n} \ar[r]^{d_n} 
&  \cdots \ar[r] & \Lambda^1P \ar[r]^{d_1:=f}\ar[d]^{\varphi_1} 
&  A \ar[r]^{d_0} \ar[d]^{\varphi_0}&  \frac{A}{IJ}\ar[r] \ar@{-->}[d]^{\overline{\varphi}_{IJ}}& 0\\
0\ar[r] &  \Lambda^0P^*\ar[r]
&  \cdots \ar[r] & (\Lambda^{n-1}P)^* \ar[r]_{(d^*)_1} &  End\left(\left|P\right|\right) \ar[r] & End\left(\frac{\left|P\right|}{IJ\left|P\right|} \right) 
  \ar[r] & 0\\   
 \enddiagram 
\end{equation}
Of course, $\left|P\right|= \Lambda^0P^*=  \Lambda^nP$. 
The rectangle at degree $r, r-1$ is given by
\begin{equation}\label{673GenRect}
 \diagram 
\Lambda^r P \ar[r]^{d_r} \ar[d]_{\varphi_r}& \Lambda^{r-1}P \ar[d]^{\varphi_{r-1}}\\
\Lambda^{n-r}P^* \ar[r]_{(d^*)_r} & \Lambda^{n-r+1}P^*\\
 \enddiagram 
\end{equation} 
We further explain the set up:
 \bE
 \item $d_r(p_1\wedge \cdots \wedge p_r)=\sum_{i=1}^r(-1)^{i-1}f(p_i) p_1\wedge \cdots \wedge p_{i-1} \wedge p_{i+1} \wedge \cdots \wedge p_r$.
 \item We follow the sign conventions in \cite {S1}, \cite[Rem. 13]{S2}, while we use homology complexes. %
 So, $(d^*)_r=(-1)^{r-1}d_{n-r+1}^*$ and the double dual identification is given by $(-1)^{r(n-r)}: \Lambda^r P \iso \Lambda^{n-r} P^*$. 
 \item The maps $\varphi_r: \Lambda^rP \lra (\Lambda^{n-r}P)^*$ are obtained by the perfect duality
 $$
 \left\{\begin{array}{l}
 \Lambda^rP \otimes \Lambda^{n-r}P \lra \left|P\right| \quad {\rm sending} \qquad \qquad \cite[Ex. 5.16(b), pp. 127]{H}\\
 (p_1\wedge \cdots  \wedge p_r)\otimes (p_{r+1}\wedge \cdots  \wedge p_n)\mapsto 
 p_1\wedge \cdots  \wedge p_r\wedge p_{r+1}\wedge \cdots  \wedge p_n\\
 \end{array}\right.
 $$
 \item The isomorphisms $\varphi_0$ and $\overline{\varphi}_{IJ}$ are given by the natural multiplication maps. 
 \item  We split $\overline{\varphi}_{IJ}=\overline{\varphi}_I\oplus \overline{\varphi}_J: \frac{A}{I}\oplus \frac{A}{J} \iso 
 End\left(\frac{\left|P\right|}{I\left|P\right|} \right)\oplus End\left(\frac{\left|P\right|}{J\left|P\right|}\right)$.  Combining with the natural identifications, we have 
 $ \f_{IJ}=\f_I\oplus \f_J$, as in the commutative diagram:
\begin{equation}\label{883fijisfiofj}
 \diagram 
  \frac{A}{I}\oplus \frac{A}{J} \ar[d]_{\overline{\varphi}_I\oplus \overline{\varphi}_J}\ar@/^/[drr]^{\f_I\oplus \f_J=\f_{IJ}}&&\\
  End\left(\frac{\left|P\right|}{I\left|P\right|} \right)\oplus End\left(\frac{\left|P\right|}{J\left|P\right|}\right)\ar[rr]_{\iota(f)}&&Ext^n\left(\frac{A}{I}, \left|P\right| \right)\oplus Ext^n\left(\frac{A}{J}, \left|P\right| \right)\\
 \enddiagram
 \end{equation}
 \eE

 
%
}
\end{remark}

%

\section{Preliminaries on Chow Witt groups}\label{BckCWSec}
We recall two different formulations of the Gersten Witt complex from \cite{BW}. Let $X$ be a regular  
scheme 
over a field $k$, with $1/2\in k$ with $\dim X=d$. Let  ${\CL}$ be an invertible sheaf on $X$ and $0\leq n\leq d$ be an integer.
We denote $X^{(n)}:=\left\{x\in X: co\dim(x)=n \right\}$. 
 For $x\in X$, let $X_x:=\spec{{\CO}_{X, x}}$
and let ${\bfD}^b_{fl}\left({\SV}(X_x, {\CL}_x \right)$ denote the bounded
derived category of complexes, of projective $X_x$-modules, with finite length homologies, and duality induced by $P\mapsto Hom(P, {\CL}_x)$.
Further, let $W^n\left({\bfD}^b_{fl}\left({\SV}(X_x, {\CL}_x \right)\right)$ denote the $n$-shifted Witt group. 
For convenience, with $X$ and ${\CL}$ being understood as above, we introduce the 
following two notations
\begin{equation}\label{701UsedHHH}
\left\{\begin{array}{l}
C^{n}\left({\bf D}_{fl}^b, W\right)=\bigoplus_{x\in X^{(n)}} W^n\left({\bfD}^b_{fl}\left({\SV}(X_x, {\CL}_x \right)\right) \\
C^{n}\left(C{\BM}, W\right)=\bigoplus_{x\in X^{(n)}} W\left(C{\BM}^n\left(X_x, {\CL}_x \right)\right)\\
\end{array}\right. 
\end{equation}
For $x\in X^{(n)}$, $W\left(C{\BM}^n\left(X_x, {\CL}_x \right)\right)$ denotes the Witt group of the
exact the category $C{\BM}^n\left(X_x \right)$ with duality $M \mapsto Ext^n\left(M, {\CL}_x\right)$.
Note, under the standing regularity hypotheses, $C{\BM}^n\left(X_x \right)$ turns out to be  the category of finite length
${\CO}_{X,x}$-modules. For brevity of notations, we would often ignore to include ${\CL}$ or ${\CL}_x$ in the notations.

From the work of    Balmer and Walter \cite{B, BW} on Witt groups 
of triangulated categories, we have the following two isomorphic descriptions of the Gersten Witt 
complex, on $X$: 
\begin{equation}\label{GerstenWitt}
\diagram 
0\ar[r] & W(X, {\CL}) \ar[r] & C^{0}\left({\bf D}_{fl}^b, W\right) \ar[r] \ar[d]^{\wr}& \cdots \ar[r] & C^{n}\left({\bf D}_{fl}^b, W\right) \ar[r]^{\qquad{\partial}^n_W} \ar[d]^{\wr}& \cdots   \\
0\ar[r] & W(X, {\CL}) \ar[r] & C^{0}\left(C{\BM}, W\right) \ar[r] & \cdots \ar[r] & C^{n}\left(C{\BM}, W\right) \ar[r]_{\qquad{d}^n_W} & \cdots  \\
\enddiagram 
\end{equation} 
terminating at the $\bigoplus_{x\in X^{(d)}}$-term, and 
where $W(X, {\CL})$ denotes the Witt group of the category of ${\SP}(X)$ with duality $P\mapsto Hom(P, {\CL})$.   We would be using these notations 
in this diagram (\ref{GerstenWitt}),
in particular, the notations  ${\partial}^n_W, {d}^n_W$ for the differentials.

%
%
%
%
 %

We proceed to recall the definition of the Chow Witt group $\widetilde{CH}^n\left(X, {\CL}\right)$ of codimension $n$ (oriented) cycles,
due to Barge and Morel \cite{BM}. 
However, 
\cite{F1} provides the most comprehensive foundation 
available on Chow Witt groups. 
Before we give the definition, we need to set up some notations, for typographical reasons. 
  For $x\in X$, denote the residue field, at $x$, by $\kappa(x)$, and $u(\kappa(x))=K_1(\kappa(x))$  denotes the group of units in $\kappa(x)$. 
  Further, for $x\in X^{(k)}$, denote 
  $$
\left\{\begin{array}{l} 
I(X_x)\subseteq W\left(C{\BM}^{k}(X_x), {\CL}_x\right),\quad{\rm the~fundamental~ideal}\\
 I^n(X_x)=I(X_x)^n  \qquad \forall n\geq 1\\
 I^n(X_x)= W\left(C{\BM}^{k}(X_x), {\CL}_x\right) \qquad \forall n\leq 0\\
\end{array}\right. 
$$
So, $\forall x\in X^{(k)}$, by \cite[Cor. 1.6]{L}, \cite[Cor. 2.3]{L}, we have 
\begin{equation}\label{kwyEqn}
\frac{I^0(X_x)}{I^1(X_x))} =\frac{{\BZ}}{2{\BZ}}, \qquad 
\frac{u(\kappa(x))}{u(\kappa(x))^2} \iso \frac{I^1\left(X_x\right)}{I^2\left(X_x\right)}
\end{equation} 
%
We ignore the coordinate associated to the duality. Fix an integer $0\leq n\leq d$. We would  define  $\widetilde{CH}^n\left(X, {\CL}\right)$,
as cohomology of a certain complex. With that in mind, 
denote the following, which would be degree $n-1, n, n+1$ terms of variety of complexes: 
\begin{equation}\label{29KIsubsup}
\left\{\begin{array}{l} 
{\SK}^{n-1}_1(X)= \bigoplus_{x\in X^{(n-1)}}  u(\kappa(x))= \bigoplus_{x\in X^{(n-1)}}  K_1(\kappa(x)) \\
{\SK}^n_0(X) =\bigoplus_{x\in X^{(n)}} K_0(\kappa(x)) =: \bigoplus_{x\in X^{(n)}} {\BZ}[x]\\
{\CI}^{n-1}_{ 1}(X)= \bigoplus_{x\in X^{(n-1)}}I^1\left(X_x\right), \quad {\CI}^{n-1}_{2}(X)= \bigoplus_{x\in X^{(n-1)}}I^2\left(X_x\right)\\
{\CI}^{n}_{0}(X)= \bigoplus_{x\in X^{(n)}}I^0\left(X_x\right), 
\quad 
{\CI}^{n}_{1}(X)= \bigoplus_{x\in X^{(n)}}I^1\left(X_x\right)\\ 
{\CI}^{n+1}_{-1}(X)=  \bigoplus_{x\in X^{(n+1)}}I^{-1}\left(X_x\right)= \bigoplus_{x\in X^{(n+1)}}W\left(X_x\right)\\
\end{array}\right. 
\end{equation} 
Thus, by (\ref{kwyEqn}), we have isomorphisms 
$$
\left\{\begin{array}{l}
\frac{{\SK}^{n-1}_1(X)}{2{\SK}^{n-1}_1(X)} =\bigoplus_{x\in X^{(n-1)}} \frac{u(\kappa(x))}{u(\kappa(x))^2}
\iso \frac{{\CI}^{n-1}_1(X)}{{\CI}^{n-1}_2(X)}\quad 
{\rm and} \\
\frac{{\SK}^{n}_0(X)}{2{\SK}^{n}_0(X)} = \bigoplus_{x\in X^{(n)}} \frac{{\BZ}}{2{\BZ}}  \iso  \frac{{\CI}^{n}_0(X)}{{\CI}^{n}_1(X)} \\
\end{array}\right.
$$
%
%
Now consider the diagram:
\begin{equation}\label{DefDia804}
\diagram
G_1^{n-1}(X) \ar[ddd]\ar[rr]^{\zeta_1} \ar@{-->}[dr]_{d^{n-1}_G}&& {\CI}^{n-1}_1(X)\ar[ddd]\ar[dr]^{d^{n-1}_I}&&\\
&G_0^n(X) \ar[ddd]\ar[rr] \ar@/_/@{-->}[drrr]_{d^n_G} && {\CI}^n_0(X)\ar[ddd]\ar[dr]^{d^n_I}&\\
&&%
&& {\CI}^{n+1}_{-1}(X)\\
{\SK}^{n-1}_1(X) 
\ar[r] \ar[dr]_{d^{n-1}_K} &\frac{{\SK}^{n-1}_1(X)}{2{\SK}^{n-1}_1(X)} 
\ar[r]_{\sim}\ar[dr] &\frac{{\CI}^{n-1}_1(X)}{{\CI}^{n-1}_2(X)}
\ar[dr]&&\\
&{\SK}^n_0(X) 
\ar[r] \ar[dr]& \frac{{\SK}^{n}_0(X)}{2{\SK}^{n}_0(X)} 
\ar[r]_{\sim}  & \frac{{\CI}^{n}_0(X)}{{\CI}^{n}_1(X)} 
&\\
&&0   &  &\\
\enddiagram
\end{equation}
%
%
%
%
%
%
%
%
%
The differential $d^{n-1}_K$ at the lower left corner is the usual order map.  
The upper right complex is, the degree $n-1, n, n+1$ portion, of a properly constructed subcomplex 
of the Gersten Witt complex (\ref{GerstenWitt}).  
%
The upper left diagonal of the diagram is obtained by cartesian product, of the lower left and upper right complex. We define three groups:
\begin{equation}\label{869CWhitt}
\left\{\begin{array}{l}
CH^n(X)=co\ker(d_K^{n-1}), \\
 C{\CI}^n\left(C, {\CL}\right)= \frac{\ker\left(d^n_I\right)}{Image\left(d^{n-1}_I\right)}, \\\ 
\widetilde{CH}^n\left(X, {\CL}\right)= \frac{\ker\left(d^n_G\right)}{Image\left(d^{n-1}_G\right)} \\
\end{array}\right.
\end{equation}
The first one is the Chow group of codimension $n$-cycles \cite{Fu}, and the last one $\widetilde{CH}^n\left(X, {\CL}\right)$ is to be called the {\bf Chow Witt groups},
({\it of co dimension $n$ oriented cycles}). 
\begin{remark}{\rm 
The following are some   information on 
the diagram (\ref{DefDia804}):
\bE
\item In   \cite{BM}, for  integers $0\leq n\leq d$,  a complex $G^{n, \bullet}$ of length $d$ was defined. 
({\it Since we fixed $n$, diagram (\ref{DefDia804}) shows only one coordinate, in the superscript  $\diagram G^{n-1}_1(X) \ar[r] & G^n_0(X) \ar[r] & \cdots \\ \enddiagram$,
indicating the degree.
The additional subscript is also helpful.})
The upper left diagonal of 
(\ref{DefDia804}), is the degree $n-1, n, n+1$ portion, of $G^{n, \bullet}$. The complex $G^{n, \bullet}$ was defined to be the cartesian product of 
the Milnor $K$-theory complex \cite{Mi}, at the lower left corner ({\it which terminates as shown in diagram \ref{DefDia804}}), 
and a subcomplex ${\CI}^{\bullet}_{\bullet -n}(X)$ of the  Gersten Witt complex (\ref{GerstenWitt}). 

However, our interest remains limited to the  degree $n-1, n, n+1$ portion of $G^{n, \bullet}$, in which case the complexes are fairly transparent and elementary. 
This is similar to the fact that, the Chow group $CH^n(X)$ is defined by the tail end of the Milnor $K$-theory complex on the lower left corner. 

%
 %



\item There are natural decompositions $G^n_0\left(X, {\CL} \right)\cong  \bigoplus_{x\in X^{(n)}} G^n_0\left(X_x\right)$,  and 
 $G^{n-1}_1\left(X, {\CL} \right)\cong  \bigoplus_{x\in X^{(n-1)}} G^{n-1}_1\left(X_x\right)$

%

\item There are other descriptions of the complex $G^{n, \bullet}\left(X, {\CL}\right)$, of the Chow Witt groups, as follows:
\bE
\item 
Given a field $F$, with $1/2\in F$, analogous to Milnor $K$-groups $K_r(F)$, 
there are groups $K_r^{MW}(F)$, known as Milnor Witt groups \cite{Mo} . Using the Milnor Witt groups $K_r^{MW}(\kappa(x))$, one can
construct Gersten complexes $C^{n, \bullet}(X,  {\CL}, K^{MW})$, exactly as in the lower right corner of the diagram (\ref{DefDia804}). One can prove that 
  $G^{n, \bullet} \left(X, {\CL}\right)\cong C^{n, \bullet}(X, K^{MW}, {\CL})$. See \cite{Mo}, \cite[Theorem 1.2, pp. 6]{F2} for further details.
 
 \item Chow-Witt group $\widetilde{CH}^n(X, {\CL})$ can also be defined as the homology of the Gersten complexes of the Grothendieck Witt groups \cite{FS,S1, M1}. 
 
 \eE
 \eE
 }
\end{remark}
 %

\vspace{3mm}
We  remark on variety of descriptions of Witt groups of regular local rings.
\begin{remark}\label{rmOnDesfWCH}{\rm 
Consider a regular local ring $(R, \m)$, with $\dim R=d$ and $1/2\in R$. Let $k=R/{\m}$. Given these, we can associate a number of Witt groups, as follows:
\bE
\item Most fundamental is the Witt group $W(k)$ of quadratic forms \cite{L}.
\item The Witt group $W({\SV}(k), *)$ of the category ${\SV}(k)$, of finite dimensional $k$-vector spaces, with duality $V \mapsto Hom(V, k)$.
\item The Witt group $W\left(C{\BM}^d(R), Ext^d\right)$ 
of the category $C{\BM}^d(R)$,
of finite length $R$-modules, with duality $M \mapsto Ext^d(M, A)$.
\item The Witt group $W\left({\bfD}^b\left(C{\BM}^d(R)\right), Ext^d\right)$ 
of the bounded derived category ${\bfD}^b\left(C{\BM}^d(R)\right)$,
of the category $C{\BM}^d(R)$, of finite length
$R$-modules, with duality induced by $M \mapsto Ext^d(M, A)$. 
\item The Witt group $W^d\left({\bfD}^d_{fl}\left({\SP}(R)\right), Hom(-, A)\right)$ of the bounded derived category ${\bfD}^d_{fl}\left({\SP}(R)\right)$,
 of complexes finite rank free $R$-modules, with finite length homology, and   duality induced by $P \mapsto Hom(P, A)$. This was mentioned above 
 (\ref{701UsedHHH}). 
({\it We would not consider skew dualities.})  
\eE
We tried to be consistent with the notations in \cite{BW}. Readers are referred to \cite{BW} for further details on derived Witt groups.
It turns out that all these Witt groups are (naturally) isomorphic \cite{B, BW, QSS, L}.
The isomorphism of the two descriptions of Gersten Witt complexes (\ref{GerstenWitt}), mentioned above, is a consequence of this fact. 
 However, some of the isomorphisms would depend on some choices to be made. 
While the first one $W(k)$ would be elementary \cite{L}, 
the last one $W^d\left({\bfD}^b_{fl}\left({\SP}(R)\right), Hom(-, A)\right)$ may also be nicer to work with. This is because 
the duality $P \mapsto Hom(P, A)$ is much more tangible than, dualities associated to "Ext". 

}
\end{remark}

 We restate \cite[Lemma 10.3.4]{F1}, as follows.

\bL\label{dx2yregsequence}{\rm 
Let $X=\spec{A}$ be a regular  scheme, $1/2\in A$,
and ${\CL}$ be an invertible sheaf on $X$. 
Consider the duality induced by $P\mapsto Hom(P, {\CL})$ on the bounded derived category ${\bfD}^b({\SP}(X))$ of the category 
${\SP}(X)$ of projective $A$-modules of finite rank.
For $x:={\wp}\in X$ denote $X_x:= \spec{A_{\wp}}$.
Let 
\begin{equation}\label{diff589}
\partial_W^n: \bigoplus_{x\in X^{(n)}} W^{n}\left({\bfD}_{fl}^{b}\left(X_x, {\CL}_x\right)\right) \lra  \bigoplus_{x\in X^{(n+1)}} W^{n+1}\left({\bfD}_{fl}^{b}\left(X_y, {\CL}_y\right)\right)
\end{equation} 
denote the differential of the Gersten Witt comples (\ref{GerstenWitt}) 
 ({\it which is a map between Witt groups of two quotient categories}).
 Let $f_1, \ldots, f_n\in A$ be a regular sequence and  $K(f_1, \ldots, f_n)$ be 
the Koszul complex. Let $\varphi: K(f_1, \ldots, f_n) \iso K(f_1, \ldots, f_n)^*$ be the corresponding symmetric form.
 Let $t\in A$  be a non zero divisor (or isomorphism)  on $\frac{A}{(f_1, \ldots, f_n)}$. 
Then, %
\begin{equation}\label{dx2yformula}
\partial^n_W: \left[\left(K(f_1, \ldots, f_n), t\varphi \right)\right]= \left[\left(K\left(f_1, \ldots, f_n; t \right), \varphi \wedge t \right)\right]
\end{equation}
where $\varphi \wedge t:K\left(f_1, \ldots, f_n; t \right) \iso K\left(f_1, \ldots, f_n; t \right)^* $ denotes the symmetric form. Evidently, $[(-, t\varphi)]$ and 
$[(-, \varphi\wedge t)]$ denote the respective elements in the Witt groups of the quotient categories.
}
\eL 
\pf 
Consider the multiplication map $t: K(f_1, \ldots, f_n) \lra K(f_1, \ldots, f_n)$.
The cone of this map  is the Koszul complex  $K\left(f_1, \ldots, f_n; t\right)$. Now, the lemma
follows from the following Lemma \ref{Descripdx2y}.
$\eop$

%

The following is a more general formulation of Lemma \ref{dx2yregsequence}. 
\bL\label{Descripdx2y}{\rm 
Let $X=\spec{A}$ be a regular  scheme, $1/2\in A$, and $\dim X=d$.
Let ${\CL}$ be an invertible sheaf on $X$. 
Consider duality induced by $P\mapsto Hom(P, {\CL})$ on  ${\bfD}^b({\SP}(X))$, 
and on the other associated categories.
Fix an integer $n$. Let $\partial^n_W$ and  $d_W^n$ denote the differentials  in the Gersten Witt complex (\ref{GerstenWitt}). 
Let ${\SD}^n({\SP}(X)) \subseteq {\bfD}^b({\SP}(X))$ denote the  subcategory of objects  $Q_{\bullet}$, 
with  $co\dim\left(H_i\left(Q_{\bullet}\right) \right) \geq n, ~\forall~i$.
In ${\SD}^n({\SP}(X))$, consider a complex $P_{\bullet}$, as follows:
$$
\begin{array}{l}
\diagram 
0\ar[r] & P_n \ar[r] & P_{n-1} \ar[r] & \cdots \ar[r] P_1 \ar[r] & P_0 \ar[r] & 0\\
\enddiagram 
\\
{\rm with~} H_i\left(P_{\bullet}\right)=0~\forall~i\neq 0,~{\rm and}~M:=H_0\left(P_{\bullet}\right)~~
{\rm Then}~~M\in C{\BM}^n(X).\\
\end{array}
$$
Let $\varphi:P_{\bullet} \iso P_{\bullet}^*$ be a symmetric (quasi) isomorphism. 
Then, $\varphi$ induces a symmetric isomorphism 
$\overline{\varphi}: M\iso M^{\vee}$  in $C{\BM}^n(X, {\CL})$ where $M^{\vee}:=Ext^n\left(M, {\CL}\right)$.

%
%
Let $t\in A$ be a non zero divisor in $M$.
The map
$$
\diagram 
P_{\bullet}\ar[d]_{\varphi t} &
0\ar[r] & P_n \ar[r]  \ar[d]^{\varphi_nt}& P_{n-1} \ar[r]  \ar[d]^{\varphi_{n-1}t}& \cdots \ar[r] P_1 \ar[r] \ar[d]^{\varphi_1t} & P_0 \ar[r] \ar[d]^{\varphi_0t}& 0\\
P_{\bullet}^* &
0\ar[r] & P_0^* \ar[r] & P_{1}^* \ar[r] & \cdots \ar[r] P_{n-1}^* \ar[r] & P_n^* \ar[r] & 0\\
\enddiagram 
$$
is a symmetric morphism in ${\SD}^{n}({\SP}(X))$, which is a lift of its image in $\frac{{\SD}^{n}({\SP}(X))}{{\SD}^{n+1}({\SP}(X))}$.
  By  \cite[Definition 5.16]{B}, \cite[Thm 2.1, 3.1]{BW}, we have
$$
\partial_W^n\left[\left(P_{\bullet}, \varphi t\right)\right]=\left[\left(Cone\left(P_{\bullet}, \varphi t\right), \psi\right)\right]
$$
where $\psi: Cone\left(P_{\bullet}, \varphi t\right) \lra Cone\left(P_{\bullet}, \varphi t\right)^*$ denotes the induced duality map.
By ($TR3$) axiom (see \cite{B}), we have a map of the triangles: 
$$
\diagram 
P_{\bullet} \ar[r]^t \ar@{=}[d]& P_{\bullet} \ar[r] \ar[d]_{\varphi}& Cone(t) \ar[r] \ar@{-->}[d]^{\mathfrak f}& TP_{\bullet}\ar@{=}[d]\\
P_{\bullet} \ar[r]_{t\varphi}  & P_{\bullet}^* \ar[r] & Cone(t\varphi) \ar[r] & TP_{\bullet}\\
\enddiagram 
$$
Since $\varphi$ is a quasi isomorphism, ${\mathfrak f}:Cone(t) \cong Cone(t\varphi)$ is also a quasi isomorphism. So, 
 $\varphi \wedge t:=\f^*\varphi\f: Cone(t) \iso Cone(t)^*$ is a symmetric isomorphism in ${\SD}^n({\SP}(X))$. Consequently, 
 ${\mathfrak f}: (Cone(t), \varphi\wedge t)\iso (Cone(\varphi t), \psi)$ is an isometry in  ${\SD}^{n+1}({\SP}(X))$. In fact, at degree $r$:
 $$
 \left\{\begin{array}{ll}
 Cone(t)_r=P_{r-1}\oplus P_r &  Cone(t)^*_r=P^*_{n-r+1}\oplus P_{n-r}^*\\ 
 Cone(\varphi t)_r=P_{r-1}\oplus P_r &  Cone(\varphi t)^*_r=P^*_{n-r+1}\oplus P_{n-r}^*\\
 \end{array}\right.
 $$
  One checks,
  $\varphi\wedge t: Cone\left(t\right) \lra \left(Cone(t)\right)^*$  has a natural description, similar to the natural dualities in Koszul complexes. 
({\it We remark that cones are very similar to Koszul complexes. For this reason the degree $r$ term of the cone can be viewed as $\left(Cone(P_{\bullet}, t)\right)_r
=P_r\oplus P_{r-1}= P_r\oplus P_{r-1}\wedge Ae_0$.})
Therefore, 
\begin{equation}\label{dWnavarphi}
\partial_{W}^n\left[\left(P_{\bullet}, t \varphi\right)\right]
= \left[\left(Cone\left(P_{\bullet}, t\varphi\right), \psi\right)\right]= \left[\left(\left(Cone(P_{\bullet}, t\right), \varphi\wedge t\right)\right]
\end{equation}
%
Further, since $t$ is a non zero divisor on $M=H_0(P_{\bullet})$, the sequence  
$$
\left\{\begin{array}{l}
\diagram 
0\ar[r] & H_0\left(P_{\bullet}\right) \ar[r]^t &H_0\left(P_{\bullet}\right) \ar[r] & H_0\left(Cone(P_{\bullet}, t)\right) \ar[r] & 0\\
\enddiagram\quad {\rm is~exact,}\\
{\rm and}\quad H_i\left(Cone(P_{\bullet}, t)\right)=0~\forall~i\neq 0\\
\end{array}\right.
$$
We obtain the commutative diagram
\begin{equation}\label{indPsi1128}
\diagram 
0\ar[r] & M \ar[r]^t \ar[d]_{\overline{\varphi}} & M \ar[r] \ar[d]^{\overline{\varphi}} & \frac{M}{tM} \ar[r] \ar[d]^{\overline{\psi}}& 0\\
0\ar[r] & Ext^n\left(M, {\CL}\right) \ar[r]_t &Ext^n\left(M, {\CL}\right) \ar[r] & Ext^{n+1}\left(\frac{M}{tM}, {\CL}\right) \ar[r] & 0\\
\enddiagram 
\end{equation}
where $\overline{\psi}$ is the induced symmetric space in $C{\BM}^{n+1}\left(X, {\CL}\right)$. Consequently,
\begin{equation}\label{notPartialWnavarphi}
d_{W}^n\left[\left(M, t \overline{\varphi}\right)\right]%
= \left[\left(\frac{M}{tM}, \overline{\psi}\right)\right]
\end{equation}

}
\eL 
\pf As elaborated! $\eop$




\section{The convergence of two obstructions} \label{ConvSec}
In this final section, we  prove our main result by establishing  a natural map $\Theta_P: \pi_0\left({\mathcal LO}(P)\right) \lra \widetilde{CH}^{n}(X)$.
\bD\label{mapDef}{\rm 
Let $X=\spec{A}$ be a Cohen-Macaulay affine scheme, with $\dim X=d$ and $1/2\in A$. Let $P$ be a projective $A$-module with $rank(P)=n$.
We adapt all the notations in above (section \ref{BckCWSec}) with ${\CL}=\Lambda^nP=\left|P\right|$, the determinant. 
 Refer to the 
diagram (\ref{DefDia804}) of definition of Chow Witt groups $\widetilde{CH}^n\left(X, \left|P\right|\right)$.
For $(I, \omega) \in \widetilde{{\mathcal LO}}(P)$,  consider 
the symmetric isomorphism $\Phi(I, \omega)=\left(\frac{A}{I}, \f_I \right)$ in $C{\BM}^n(X, \left|P\right|)$, as defined in  
Definition \ref{74DefPhira}, item (\ref{810ItemPhiIomeg}), and equation (\ref{801MoAss}). Recall, that the same notation $\Phi(I, \omega)$ was used to denote the symmetric form and 
its isometry class in the monoid ${\bf M}C{\BM}^n(X, \left|P\right|)$.
This defines an elements 
$$
\left\{\begin{array}{l}
\Phi(I, \omega)_I \in {\CI}^n_0(X) =\bigoplus_{x\in X^{(n)}}I^0(X_x)\\
\Phi(I, \omega)_K \in  {\SK}^n_0(X)= \bigoplus_{x\in X^{(n)}}K_0(X_x)\\
\end{array}\right. 
$$  
They patch to define an element in $\Phi(I, \omega)_G \in G^n_0\left(X, \left|P\right|\right)$.
Note $\Phi(I, \omega)_I$ comes from a symmetric space in $C{\BM}^n(X)$, or equivalently from ${\SD}^n({\SP}(X))$.
Consequently, $d_I^n\left(\Phi(I, \omega)_I\right)=0$, 
which follows from construction of the Gersten Witt complex  (\ref{GerstenWitt},  \cite{BW}). Hence   $d^n_G\left(\Phi(I, \omega)_G\right)=0$.
%
Therefore, $\Phi(I, \omega)_G$ represents an element in $\widetilde{CH}^n\left(X, \left|P\right|\right)$.
Define, a set theoretic map 
\begin{equation}\label{OmegaMap}
\Omega: \widetilde{{\mathcal LO}}(P) \lra \widetilde{CH}^n\left(X, \left|P\right|\right) ~~ {\rm by}~~
\Omega(I, \omega):=\overline{\Phi(I, \omega)_G} \in  \widetilde{CH}^n\left(X, \left|P\right|\right)
%
\end{equation}
The following commutative diagram would be helpful,
\begin{equation}\label{DefOmegaPh519}
\diagram 
\widetilde{{\mathcal LO}}(P) \ar[r]^{\Phi\qquad} \ar@/_/[drr]_{\Omega}& {\bf M}\left(C{\BM}^n\left(X, \left|P\right|\right)\right)\ar[r]^{\qquad \iota_n} 
& \ker(d^n_G) 
\ar[d]^{q}\\
&& \widetilde{CH}^n\left(X, \left|P\right|\right)\\
\enddiagram
\end{equation}
where 
$\iota_n$ is a the map obtained by point wise localization map. 
}
\eD

\vspace{3mm}
\begin{remark}{\rm 
Chow Witt groups $\widetilde{CH}^n\left(X, \left|P\right|\right)$ can also be defined, equivalently, as the cohomology of the Gersten complex of the Grothedieck Witt groups.
So,  there is a direct $GW$-way to look at all these  \cite[Remark 4.14(1), 4.5(1)]{M1}. 
}
\end{remark}

\vspace{3mm}
We state the following on homotopy invariance of Chow Witt groups.
\bP\label{homoInv1Jan21}{\rm 
Let $X$ be a regular  scheme over a field $k$, with $1/2\in k$, $d=\dim X$. Fix an integer $0\leq n\leq d$. 
and  an invertible sheaf ${\CL}$ on $X$. All the dualities considered are induced by the duality $Q\mapsto Hom(Q, {\CL})$ on ${\SP}(X)$. 
Let $p:X\times {\BA}^1 \lra X$ be the projection map. Then, $p$ induces an isomorphism 
$$
p^*:\widetilde{CH}^n(X, {\CL}) \iso \widetilde{CH}^n(X\times {\BA}^1, p^*{\CL})\qquad {\rm of ~the~Chow~ Witt ~groups.} 
$$
}
\eP 
\pf See \cite[Cor. 11.3.3]{F1}. $\eop$

\vspace{3mm}
The following is the main theorem in this article. 
\bT\label{needsProof}{\rm 
Suppose $A$ is a regular ring, containing a field $k$  with $1/2\in k$, and $\dim A=d\geq 2$. 
Let $P$ be a projective $A$-module of rank $n$. 
Then, the map $\Omega$ in (\ref{OmegaMap}) factors through a set theoretic map $\Theta_P$, as in the 
commutative diagram:
\begin{equation}\label{theFinalMap} 
\diagram 
 \widetilde{{\mathcal LO}}(P)  \ar@/_/[dr]_{\Omega} \ar[r]^{\beta}& \pi_0\left({\mathcal LO}(P)\right)\ar@{-->}[d]^{\Theta_P}\\ 
 & \widetilde{CH}^n\left(X, \left|P\right|\right)\\
\enddiagram 
\end{equation}
}
\eT
\pf 
By  (\ref{transIso}),  $\pi_0\left({\mathcal LO}(P)\right)=\pi_0\left( \widetilde{{\mathcal LO}}(P) \right)$. 
Let $(I_0, \omega_0), (I_1, \omega_1)\in  \widetilde{{\mathcal LO}}(P)$ be such that
 $\beta((I_0, \omega_0))=\beta((I_1, \omega_1))\in \pi_0\left( \widetilde{{\mathcal LO}}(P) \right)$. 
By 
 (\ref{moviongRem}), 
there is $(I, \omega)\in {\mathcal LO}^n(P[T])$ such that $(I, \omega)_{|T=0}= (I_0, \omega_0)$ and $(I, \omega)_{|T=1}= (I_1, \omega_1)$.
The  projection map $p: X\times {\BA}^1 \lra X$ induces a pull back map 
$$
p^*:\widetilde{CH}^n(X, \left|P\right|) \iso \widetilde{CH}^n(X\times {\BA}^1, p^*\left|P\right|)$$
 which is an 
isomorphism 
(\ref{homoInv1Jan21}) or \cite[Theorem 2.15]{F2}. 
Consider the commutative diagram ({\it while we duplicate notations $\Omega$, $\Phi$ etc.})
\begin{equation}\label{LOCHComm}
\diagram 
&& \widetilde{CH}^n(X, \left|P\right|) \ar[ddd]^{p^*}_{\wr}\\
\widetilde{{\mathcal LO}}(P) \ar@/^/[rru]^{\Omega}\ar[r]_{\Phi\qquad} \ar[d]&{\bf M}\left(C{\BM}^n\left(X, \left|P\right|\right)\right)\ar[ru]_{q\iota_n}\ar[d]^{p^*}&\\ 
\widetilde{{\mathcal LO}}(P[T]) \ar@/_/[rrd]_{\Omega\qquad} \ar[r]^{\Phi\qquad}&{\bf M}\left(C{\BM}^n\left(X\times {\BA}^1, p^*\left|P\right|\right)\right)\ar[rd]_{q\iota_n} & \\ 
& & \widetilde{CH}^n(X\times {\BA}^1, p^*\left|P\right|) \\
\enddiagram
\end{equation}
%
%
We are required to prove that $\Omega(I_0, \omega_0)= \Omega(I_1, \omega_1)$. Since $p^*$ (the last vertical arrow) is an isomorphism, it is enough to prove that 
\begin{equation}\label{homoSaman}
p^*\left(\Omega(I_i, \omega_i)\right))= \Omega(I, \omega) \qquad \qquad {\rm for}\quad i=0, 1
\end{equation}
Since $T$ and $T-1$ are interchangeable, it would suffice to establish the case $i=0$. This is established below,   in Lemma \ref{contraction550}.
\pic $\eop$


\subsection{The retraction  $T=0$}
We establish the equation (\ref{homoSaman}), in this subsection.
\bL\label{contraction550}{\rm 
Let $X=\spec{A}$ be a regular affine scheme, over $\spec{k}$, where $k$ is a field with $1/2\in k$ and  $\dim X=d$.
Let $P$ be a projective $A$-module with $rank(P)=n$, with $2\leq n \leq d$. 
Let $(I, \omega)\in {\mathcal LO}^n(P[T])$ be such that $(I_0, \omega_0):=(I, \omega)_{|T=0} \in \widetilde{{\mathcal LO}}(P)$. 
Then, $p^*\left(\Omega(I_0, \omega_0)\right)= \Omega(I, \omega)$, with notations as in (\ref{LOCHComm}). 
}
\eL


%

\pf 
Let $f_0:P \sur I_0\cap J_0$ be a lift of $\omega_0:P \sur \frac{I_0}{I_0^2}$ such that $height(J_0)\geq n$, and $I_0+J_0=A$. 
 By \cite[Lem. 2.2]{M2}, there is a lift  $h(T):P[T] \lra I$ of $\omega$ such that $h(0)=f_0$. Now, $I=(h(P[T]), s)$, with $s\in I^2$. 
 Consider the set,
 $V=\{\wp\in \spec{A[T]}: height(\wp)\leq n-1, sT\notin {\wp} \}$. Now, $(h, sT)\in P[T]^*\oplus A[T]$ is basic on $V$. By basic element 
 theory  \cite{M3}, we have $f:=h+sTg$ is basic on $V$, for some $g\in P[T]^*$ . It follows, (1) $f(0)=f_0$ and (2) $f$ is a lift of $\omega$ and (3) $height(f(P[T]))\geq n$. So,
 we can write $f(P[T])=I\cap J$, with $height(J)\geq n$ and $I+J=A[T]$. We fix such a lift $f$ of $\omega$. Also, let $\omega_J:P[T] \sur \frac{J}{J^2}$, 
 $\omega_{IJ}:P[T] \sur \frac{IJ}{(IJ)^2}$, 
  be induced by $f$. 
 Then, $(J, \omega_J), (IJ, \omega_{IJ}) \in \widetilde{{\mathcal LO}}(P[T])$.
 Likewise, denote $(J_0, \omega_{J_0}), (I_0J_0, \omega_{I_0J_0}) \in \widetilde{{\mathcal LO}}(P)$.
The Koszul complex of $f$ anf $f_0$ leads to the  symmetric spaces,  which may be considered in the $C{\BM}^n(-)$-categories or in the derived categories:
$$
\diagram 
0\ar[r] & \Lambda^0P[T]^*\ar[r] &\cdots \ar[r] & \Lambda^{n-1}P[T]^* \ar[r] & \Lambda^nP[T]^* \ar[r] &  Ext^k\left(\frac{A[T]}{(I\cap J)}, {\CL}[T]\right) \ar[r] & 0\\
0\ar[r] & \Lambda^nP[T] \ar[r] \ar[u]_{\varphi}&\cdots \ar[r]  
& \Lambda^1P[T] \ar[r]_f \ar[u]_{\varphi} & \Lambda^0P[T] \ar[r] \ar[u]_{\varphi}&   \frac{A[T]}{(I\cap J)} \ar[r]\ar[u]_{\f_I\oplus \f_J} & 0\\
0\ar[r] & \Lambda^nP \ar[r] \ar[d]^{\varphi(0)}&\cdots \ar[r]  
& \Lambda^1P \ar[r]^{f_0} \ar[d]^{\varphi(0)} & \Lambda^0P \ar[r] \ar[d]^{\varphi(0)} &   \frac{A}{(I(0)\cap J(0))} \ar[r]\ar[d]^{{\f}_{I_0}\oplus {\f}_{J_0}} & 0\\
0\ar[r] & \Lambda^0P^*\ar[r] &\cdots \ar[r] & \Lambda^{n-1}P^* \ar[r] & \Lambda^nP^* \ar[r] &  Ext^k\left(\frac{A}{(I(0)\cap J(0))}, {\CL}\right) \ar[r] & 0\\
\enddiagram 
$$
Denote $X[T]:=X\times {\BA}^1=\spec{A[T]}$.  
We proceed to use Lemma \ref{Descripdx2y}, (\ref{indPsi1128}, \ref{notPartialWnavarphi}). 
 Denote the symmetric spaces in ${\SD}^{n+1}({\SP}(X[T])$, as follows
$$
\left\{\begin{array}{ll}
{\bfC}(T\varphi):=(Cone(T), \varphi\wedge T) \cong Cone(T\varphi) &{\rm cone~of}~ T\varphi:\Lambda^{\bullet}P[T] \lra \Lambda^{\bullet}P[T]\\ 
{\bfC}(T\varphi(0)[T]):=(Cone(T), \varphi(0)[T]\wedge T)
& {\rm cone~of}~ T\varphi(0)[T]:\Lambda^{\bullet}P[T] \lra \Lambda^{\bullet}P[T]\\
\end{array}\right.
$$
We claim that these  symmetric spaces ${\bfC}(T\varphi)$, ${\bfC}(T\varphi(0)[T])$ are isometric. Write $f(T)=f_0+Tf_1+\cdots+T^mf_m$, with $f_i\in Hom(P, A)$. 
Consider the map 
$$
\Delta=\left(\begin{array}{cc}1_{A[T]} &   \frac{f-f_0[T]}{T}\\ 0 & 1_{P[T]}\\ \end{array}\right): A[T]\oplus P[T] \iso A[T]\oplus P[T] 
$$
$$
{\rm Then,} \qquad \qquad 
\left(\begin{array}{cc} -T & f \\ \end{array}\right)
\left(\begin{array}{cc}1_{A[T]} &   \frac{f-f_0[T]}{T}\\ 0 & 1_{P[T]}\\ \end{array}\right)
=
\left(\begin{array}{cc} -T & f_0[T] \\ \end{array}\right)
$$
Therefore, the diagram
$$
\diagram 
A[T]\oplus P[T] \ar[rr]^{\qquad(-T,  f_0[T])}\ar[d]_{\Delta} && A[T] \ar[r]\ar[d]^1 & \frac{A[T]}{IJ} \ar[r]\ar[d]^1 & 0\\
A[T]\oplus P[T] \ar[rr]_{\quad(-T, f)} && A[T] \ar[r] & \frac{A[T]}{IJ} \ar[r] & 0\\
\enddiagram
\qquad {\rm commutes.} 
$$
We write $Q=A[T]\oplus P[T]$ and $Q_r=\Lambda^r Q$. Then, the cone ${\bfC}(T\varphi)$ is the Koszul complex of the map $(-T, f):Q\lra A[T]$,
and ${\bfC}(T\varphi(0)[T])$ is the Koszul complex of the map $(-T, f(0)[T]):Q\lra A[T]$. 

$$
{\rm Write}\qquad
\left\{\begin{array}{l}
M =\frac{A[T]}{(IJ, T)}=
\frac{B}{(I_0 J_0, T)}\\
M^{\vee}=Ext^{n+1}\left(M, \left|P[T]\right|\right)\\
\end{array}\right.
$$
Consider the following composition of maps of Koszul complexes, and duality:
$$
\diagram 
0 \ar[r] &\Lambda^{n+1} Q \ar[r]^{\partial_0}\ar@{=}[d] &
\cdots \ar[r] & \Lambda^2 Q \ar[r] \ar[d]^{\Lambda^2\Delta}
& Q \ar[r]^{(-T, f_0[T])}\ar[d]^{\Delta} & A[T] \ar[r]\ar@{=}[d] & M \ar@{=}[d]\ar[r] & 0\\
0 \ar[r] &\Lambda^{n+1} Q \ar[r]_{\partial_1}\ \ar[d]^{can}&
\cdots \ar[r] & \Lambda^2 Q \ar[r]\ar[d]^{can} & Q \ar[r]_{(-T, f)} \ar[d]^{can} & A[T] \ar[r] \ar[d]^{can=1}& M \ar[r]\ar[d]^{\f_I\oplus {\f}_J} & 0\\
0 \ar[r] &  A[T]^*\ar[r]\ar@{=}[d] &
\cdots \ar[r] & \Lambda^{n-1} Q^* \ar[r] \ar[d]
& \Lambda^{n} Q^* \ar[r]^{\partial_1^{*}} \ar[d]^{\Lambda^{n}\Delta^*} & \Lambda^{n+1} Q^* \ar[r]\ar@{=}[d] & M^{\vee} \ar@{=}[d]\ar[r] & 0\\
0 \ar[r] &  A[T]^*\ar[r]  & 
\cdots \ar[r] & \Lambda^{n-1} Q^* \ar[r]  
& \Lambda^{n} Q^* \ar[r]_{\partial_0^{*}}  & \Lambda^{n+1} Q^* \ar[r] & M^{\vee} \ar[r] & 0\\
\enddiagram 
$$
Since $\det(\Delta)=1$,
it follows from (\ref{DetFitthBll}), the composition is the canonical map $can: \Lambda^{\bullet}Q \iso  \Lambda^{\bullet}Q^*$.
Therefore, $\Lambda^{\bullet}\Delta: {\bfC}(T\varphi(0)[T]) \iso {\bfC}(T\varphi)$  is an isometry, in ${\SD}^{n+1}({\SP}(X))$. 
%
At degree zero, the compositions agree. It follows that from this that 
$$
{\f}_I\oplus {\f}_J={\f}_{I_0}[T]\oplus {\f}_{J_0}[T]: M \iso M^{\vee}
$$
 So, $\forall ~y\in X[T]^{(n+1)}$, the images  agree:
\begin{equation}\label{670TdudiMale}
 \left(\frac{A[T]}{I}, \f_I \right)_y= \left(\frac{A[T]}{I_0[T]}, \f_{I_0}[T] \right)_y\in W\left(C{\BM}^{n+1}(X[T]_y, \left|P\right|[T]_y \right)
\end{equation}
%
%
Up to isometry, we have
$$
\left\{\begin{array}{l}
\Phi(I, \omega)=\left(\frac{A[T]}{I}, \f_I \right)\in {\bf M}\left(C{\BM}^n(X[T])\right),\\
 \Phi(I_0, \omega_0)= \left(\frac{A}{I_0}, \f_{I_0} \right)\in {\bf M}\left(C{\BM}^n(X)\right)\\
p^*  \Phi(I_0, \omega_0)=\left(\frac{A[T]}{I_0[T]}, \f_{I_0}[T] \right)\in {\bf M}\left(C{\BM}^n(X[T])\right)\\
 \end{array}\right. 
$$
We obtain $\Phi(I, \omega)_I\in {\CI}^n_0(X[T])$ and $\Phi(I_0, \omega_0)_I\in {\CI}^n_0(X)$, by point wise localization. Now,
$$
\left\{\begin{array}{l}
T\cdot \Phi(I, \omega):=\left(\frac{A[T]}{I}, T\f_I \right)\\
T\cdot p^*  \Phi(I_0, \omega_0):=\left(\frac{A[T]}{I_0[T]}, T\f_{I_0}[T] \right)\\
 \end{array}\right. 
$$
represent symmetric forms in $C{\BM}^n(X[T])$. We use the  notation 
$\langle T\rangle \cdot \Phi(I, \omega)_I$ and $\langle T\rangle \cdot p^*\Phi(I_0, \omega_0)_I$ to denote the corresponding elements in ${\CI}^n_0(X[T])$. 
Refer to  the notations in (\ref{GerstenWitt}, \ref{29KIsubsup}, \ref{DefDia804}). By Lemma \ref{Descripdx2y}, we have  
$$\left\{\begin{array}{l}
\partial_W^n\left(\Lambda^{\bullet}P[T], T\varphi\right) =\left[\left(\left(\Lambda^{\bullet}P[T], Cone(T)\right), \varphi\wedge T\right)\right]. \quad 
{\rm And,~equivalently,}\\
d^n_W\left(T\cdot\left( \Phi\left(I, \omega \right)_I\perp \Phi\left(J, \omega_J \right)_I\right)\right)=\sum_{x\in X[T]^{(n+1)}} \left[\left(\frac{A[T]}{(I\cap J, T)}, \psi \right)_x \right]
\end{array}\right.
$$
where $\psi:\frac{A[T]}{(I\cap J, T)} \iso Ext^{n+1}\left(\frac{A[T]}{(I\cap J, T)}, p^*{\CL}\right)$ is the symmetric isomorphism induced by $\varphi\wedge T$ (see \ref{indPsi1128}, \ref{notPartialWnavarphi}).
Since $I+J=A[T]$ and $I_0+J_0=A$, it follows
\begin{equation}\label{d1IPhiIOmega}
\left\{\begin{array}{l}
d^n_W\left(\langle T\rangle \cdot \Phi\left(I, \omega \right)_I\right) = \sum_{x\in X[T]^{(n+1)}} \left[\left(\frac{A[T]}{(I, T)}, \psi \right)_x \right] \\
=d^n_I\left((\langle T \rangle - \langle 1 \rangle)\cdot \Phi\left(I, \omega \right)_I\right)\\ 
\end{array}\right.
\end{equation}
in ${\CI}^{n+1}_0(X[T], p^*{\CL})$. 
Similarly, 
\begin{equation}\label{d1pStarIPhiIOm0}
\left\{\begin{array}{l}
d^n_W\left(\langle T \rangle\cdot p^*\Phi\left(I_0, \omega_0 \right)_I\right) = \sum_{x\in X[T]^{(n+1)}} \left[\left(\frac{A[T]}{(I_0, T)}, \psi_0 \right)_x \right] \\
=d^n_I\left((\langle T \rangle - \langle 1 \rangle)\cdot p^*\Phi\left(I_0, \omega_0 \right)_I\right) \\
\end{array}\right.
\end{equation}
where $\psi_0$ is induced by $\varphi(0)[T] \wedge T$ (see \ref{indPsi1128}, \ref{notPartialWnavarphi}).
It follows from (\ref{670TdudiMale}) that 
\begin{equation}\label{claim1199}
 d^n_I\left((\langle T \rangle - \langle 1 \rangle)\cdot \Phi\left(I, \omega \right)_I\right)= d^n_I\left((\langle T \rangle - \langle 1 \rangle)\cdot p^*\Phi\left(I(0), \omega_0 \right)_I\right) 
\end{equation}
in $ {\CI}^{n+1}_0(X[T], \left|P[T]\right|)$. Similarly, we have 
$$
d^n_K\left(T\cdot \Phi\left(I, \omega \right)_K\right)=d^n_K\left(T \cdot p^*\Phi\left(I(0), \omega \right)_K\right) \in {\SK}_0^{n+1}(X[T])
$$
Patching these two identities, it follows
\begin{equation}\label{754finGIdentity}
d^n_G\left(T\cdot \Phi\left(I, \omega \right)_G\right)=d^n_G\left(T \cdot p^*\Phi\left(I(0), \omega_0 \right)_G\right)
\quad \in   G^{n+1}_0\left(X[T], \left|P[T]\right| \right)
\end{equation}
Consider the homotopy invariance,
$$
p^*: \widetilde{CH}^n\left(X, \left|P\right| \right) \iso  \widetilde{CH}^n\left(X\times {\BA}^1,\left|P[T]\right| \right)
$$
It follows from (\ref{754finGIdentity}), and the definition of the retraction, that 
$$
(p^*)^{-1}\Omega(I, \omega)= (p^*)^{-1}(p^*\Omega(I_0, \omega_0)= \Omega(I_0, \omega_0)
$$
\pic $\eop$ 



\bC\label{ChCW}{\rm 
With the notations as above (\ref{needsProof}), we  have
$$
\Theta_P(\varepsilon_H(P))= \varepsilon_{CW}(P) \in \widetilde{CH}^n(X)
$$
where $\varepsilon_H(P)\in \pi_0\left(\widetilde{{\mathcal LO}}(P) \right)$ is the homotopy obstruction defined in (\ref{neutrlObsCls}),
 and $\varepsilon_{CW}(P)$ denotes the oriented Chern class, as defined in \cite{BM, F1}. 
}
\eC 
%

\vspace{3mm}
The following lemma is standard that we used above. 
\bL\label{DetFitthBll}{\rm 
Suppose $A$ is a commutative noetherian ring and $Q$ is a projective $A$-module with $rank(Q)=n$  and let ${\CL}:=\Lambda^nQ=\left|Q\right|$. 
For  $A$-modules $M$, we denote $M^*:=Hom(M, {\CL})$. 
Let $0\leq r\leq n$ and $can:\Lambda^r Q \lra \Lambda^{n-r} Q^*$ be the  isomorphism induced by the perfect duality $\Lambda^r Q\otimes \Lambda^{n-r} Q
\lra {\CL}$. Let $\Delta: Q \iso Q$ be an isomorphism, with $\det(\Delta)=1$. Then, the diagram
\begin{equation}\label{503Can33}
\diagram 
\Lambda^rQ \ar[rr]^{\Lambda^r\Delta}\ar@/_/[drrr]_{can} & &\Lambda^rQ \ar[r]^{can} & \Lambda^{n-r}Q^*\ar[d]^{ \Lambda^{n-r}\Delta^*}\\
&&&  \Lambda^{n-r}Q^*\\
\enddiagram 
\end{equation}
commutes. 
}
\eL 
\pf First, we prove this for $r=1$. Let $q\in Q$. Let $can \Delta(q)\in \Lambda^{n-1}Q^*$ is given by 
$can \Delta(q)(q_2\wedge \cdots \wedge q_n)=\Delta(q)\wedge q_2\wedge \cdots \wedge q_n$.
Then, $\lambda_q:= \Lambda^{n-r}\Delta^*can \Delta(q)$ is given by the commutative diagram 
$$
\diagram 
\Lambda^{n-1}Q \ar[rr]^{\Lambda^{n-1}Q}\ar@/_/[drr]_{\lambda_q} && \Lambda^{n-1}Q \ar[d]^{can \Delta(q)}\\
&& \Lambda^nQ\\
\enddiagram 
$$
So,
$$
\lambda_q(q_2\wedge \cdots \wedge q_n)= \Delta(q) \wedge\Lambda^{n-1} \Delta (q_2\wedge \cdots \wedge q_n)=\Lambda^n(\Delta)(q\wedge q_2\wedge \cdots \wedge q_n)
$$
Since $\det\Delta=1$, it follows that, the upper right composition (\ref{503Can33}), sends $q\mapsto (q_2\wedge \cdots \wedge q_n \mapsto q\wedge q_2\wedge \cdots \wedge q_n)$,
which is the map $can:Q \lra \Lambda^{n-1}Q^*$. 

For the general case, for $q_1, \ldots, q_r\in Q$, let $ev_{q_1q_2\ldots q_r}\in \Lambda^{n-r}Q^*$, denote the map $q_{r+1}\wedge \cdots \wedge q_n 
\mapsto q_1\wedge \cdots \wedge q_{r+1}\wedge \cdots \wedge q_n$. So, $can \Lambda^r \Delta(q)= ev_{{\Delta(q_1)\Delta(q_2)\ldots \Delta(q_r)}}$.
So, $\Lambda^{n-r}\Delta can \Lambda^r \Delta(q)$ is given by the commutative diagram
$$
\diagram 
\Lambda^{n-r}Q \ar[rr]^{\Lambda^{n-r}\Delta}\ar@/_/[rrd] && \Lambda^{n-r}Q\ar[d]^{ev_{{\Delta(q_1)\Delta(q_2)\ldots \Delta(q_r)}}}\\ 
&& \Lambda^n Q\\ 
\enddiagram
$$
which sends 
$$
q_{r+1}\wedge \cdots \wedge q_n \mapsto \Delta(q_1)\wedge \cdots \wedge\Delta(q_r)\wedge \Delta(q_{r+1})\wedge \cdots \wedge \Delta(q_n)=\Lambda^n\Delta(q_1\wedge \cdots 
\wedge q_r \wedge q_{r+1} \wedge \cdots \wedge q_n)
$$
This is precisely the map $can:\Lambda^rQ \lra \Lambda^{n-r}Q^*$. \pic $\eop$





\end{document}